\documentclass[review,onefignum,onetabnum]{siamart190516}

\usepackage{enumitem}
\usepackage{amsfonts}
\usepackage{graphicx}
\usepackage{epstopdf}
\usepackage{algorithmic}
\ifpdf
  \DeclareGraphicsExtensions{.eps,.pdf,.png,.jpg}
\else
  \DeclareGraphicsExtensions{.eps}
\fi

\numberwithin{theorem}{section}

\usepackage{hyperref}
\hypersetup{
     colorlinks   = true,
     citecolor    = blue,
     linkcolor    = blue,
     urlcolor     = blue
}
\usepackage{tikz,pgfplots}
\usetikzlibrary{shapes}
\usepackage{graphicx}
\usepackage{algorithm,algorithmic}
\newcommand{\rottext}[1]{\rotatebox{90}{\hbox to 20mm{\hss #1\hss}}}
\newcommand\norm[1]{\left\lVert#1\right\rVert}
\usepackage{datetime}

\newcommand{\R}{\mathbb{R}}

\renewcommand{\norm}[2][]{\left\Vert#2\right\Vert_{#1}}

\newcommand{\kron}{\otimes}
\DeclareMathOperator*{\argmin}{arg\,min}                  
\renewcommand{\t} {^{\top}}

\renewcommand{\phi}{\mathbf{\varphi}}

\newcommand{\bfzero}{{\bf0}}

\newcommand{\bfd}{\mathbf{d}}

\newcommand{\bfF}{\mathbf{F}}
\newcommand{\bfT}{\mathbf{T}}
\newcommand{\bff}{\mathbf{f}}
\newcommand{\bfA}{\mathbf{A}}
\newcommand{\bfC}{\mathbf{C}}
\newcommand{\bfE}{\mathbf{E}}
\newcommand{\bfU}{\mathbf{U}}

\newcommand{\bfR}{\mathbf{R}}
\newcommand{\bfI}{\mathbf{I}}
\newcommand{\bfD}{\mathbf{D}}
\newcommand{\bfn}{\mathbf{n}}

\newcommand{\bfc}{\mathbf{c}}
\newcommand{\bfH}{\mathbf{H}}

\newcommand{\bfx}{\mathbf{x}}

\newcommand{\bfe}{\mathbf{e}}
\newcommand{\bfu}{\mathbf{u}}
\newcommand{\bfy}{\mathbf{y}}

\newcommand{\bfM}{\mathbf{M}}
\newcommand{\bfG}{\mathbf{G}}
\newcommand{\bfL}{\mathbf{L}}
\newcommand{\bfw}{\mathbf{w}}
\newcommand{\bfW}{\mathbf{W}}
\newcommand{\bfz}{\mathbf{z}}

\newcommand{\bfg}{\mathbf{g}}
\newcommand{\bfh}{\mathbf{h}}
\newcommand{\bfr}{\mathbf{r}}
\newcommand{\bfs}{\mathbf{s}}
\newcommand{\bfv}{\mathbf{v}}
\newcommand{\bfV}{\mathbf{V}}

\newcommand{\bfQ}{\mathbf{Q}}

\newcommand{\bfZ}{\mathbf{Z}}

\newcommand{\calL}{\mathcal{L}}

\newcommand{\calN}{\mathcal{N}}
\newcommand{\calO}{\mathcal{O}}

\newcommand{\calR}{\mathcal{R}}

\newcommand{\bfdelta}{{\boldsymbol{\delta}}}
\newcommand{\bfepsilon}{{\boldsymbol{\epsilon}}}

\newcommand{\bfzeta}{{\boldsymbol{\zeta}}}

\newcommand{\bfmu}{{\boldsymbol{\mu}}}

\newcommand{\bfxi}{{\boldsymbol{\xi}}}

\newcommand{\bbR}{\mathbb{R}}

\usetikzlibrary{arrows}

\newlength\iwidth
\newlength\iheight

\usepackage{subcaption}
\usepackage{caption}
\usepackage{tabularx}
\usepackage{graphicx}
\usepackage{array}
\usepackage{longtable}
\usepackage{verbatim}

\newdimen\iwidth
\newdimen\iheight

\theoremstyle{plain}

\newcommand{\TheTitle}{Hybrid Projection Methods for Solution Decomposition in Large-scale Bayesian Inverse Problems}

\newcommand{\TheAuthors}{Chung, Jiang, Miller, and Saibaba}

\headers{Hybrid projection methods for solution decomposition}{\TheAuthors}

\title{{\TheTitle}
\thanks{Updated \today.}
\funding{This work was partially supported by the National Science Foundation ATD program under grants DMS-2026841, 2026830, and 2026835.  This material was also based upon work partially supported by the National Science Foundation under Grant DMS-1638521 to the Statistical and Applied Mathematical Sciences Institute. Any opinions, findings, and
conclusions or recommendations expressed in this material are those of the author(s) and do not necessarily reflect the views of the National Science Foundation.}}

\author{
  Julianne Chung\thanks{Department of Mathematics, Computational Modeling and Data Analytics Division, Academy of Integrated Science, Virginia Tech, Blacksburg, VA, USA
    (\email{jmchung@vt.edu}, \url{http://www.math.vt.edu/people/jmchung/}).}
  \and
Jiahua Jiang\thanks{Department of Mathematics, University of Birmingham, UK
  (\email{j.jiang.3@bham.ac.uk}).}
  \and
Scot M. Miller\thanks{Department of Environmental Health and Engineering, Johns Hopkins University
  (\email{smill191@jhu.edu}).}
    \and
Arvind K. Saibaba\thanks{Department of Mathematics, North Carolina State University
  (\email{asaibab@ncsu.edu}).}
}

\usepackage{amsopn}

\definecolor{darkcyan}{rgb}{0.0, 0.55, 0.55}

\usepackage[normalem]{ulem}

\begin{document}
\nolinenumbers
\maketitle

\begin{abstract}
We develop hybrid projection methods for computing solutions to large-scale inverse problems, where the solution represents a sum of  different stochastic components. Such scenarios arise in many imaging applications (e.g., anomaly detection in atmospheric emissions tomography) where the reconstructed solution can be represented as a combination of two or more components and each component contains different smoothness or stochastic properties.  In a deterministic inversion or inverse modeling framework, these assumptions correspond to different regularization terms for each solution in the sum. Although various prior assumptions can be included in our framework, we focus on the scenario where the solution is a sum of a sparse solution and a smooth solution. For computing solution estimates, we develop hybrid projection methods for solution decomposition that are based on a combined flexible and generalized Golub-Kahan processes. This approach integrates techniques from the generalized Golub-Kahan bidiagonalization and the flexible Krylov methods. The benefits of the proposed methods are that the decomposition of the solution can be done iteratively, and the regularization terms and regularization parameters are adaptively chosen at each iteration. 
Numerical results from photoacoustic tomography and atmospheric inverse modeling demonstrate the potential for these methods to be used for anomaly detection. 
\end{abstract}


\thispagestyle{plain}

\section{Introduction}
\label{sec:introduction}
In many inverse problems, the ability to efficiently and accurately detect anomalies from observed data can have significant benefits. For example, in atmospheric inverse modeling, large-scale anomalous emissions of greenhouse gasses and air pollution pose threats to human health, state emissions targets, and energy security.  We need inverse models that can identify anomalous emissions events quickly -- so the leak or malfunction in question can be fixed. 
However, inverse models that can identify anomalous emissions events require \emph{more complicated prior models}, in particular, models that can incorporate multiple complex sources with different smoothness properties. It is desirable to capture both anomalies (e.g., sparsely distributed events representing anomalous emissions like natural gas blowouts or point sources like large power plants) and smooth regions (e.g., representing broad scale emissions patterns from area sources). Moreover, the quality of the reconstruction depends crucially on the \emph{choice of appropriate hyperparameters} that govern the prior and noise distributions, and estimating these parameters prior to inversion can be prohibitively expensive.
For these and other inverse problems where the solution must capture different stochastic properties, we describe efficient and flexible iterative methods for reconstructing
solutions that have a combination of smooth regions with sparse anomalies (e.g., for anomaly detection).

More specifically, we consider linear inverse problems of the form,
\begin{equation}
 \label{eq:sdproblem}
 \bfd = \bfA\bfs + \bfdelta, \quad \mbox{with} \quad \bfs = \bfs_1+\bfs_2,
\end{equation}
where the goal is to reconstruct the desired parameters $\bfs_1, \bfs_2 \in \bbR^n$, given forward operator (or parameter-to-observable map) $\bfA \in \bbR^{m\times n}$ and the observed data $\bfd \in \bbR^m$. We assume that the measurement errors $\bfdelta$ are realizations of Gaussian random variables, i.e., $\bfdelta \sim \calN(\bfzero, \bfR)$, where $\bfR$ is a symmetric positive definite (SPD) matrix, and that $\bfs_1$ and $\bfs_2$ are mutually independent and are realizations from \textit{different} distributions. Contrary to most inverse problems that involve estimating the unknown parameters $\bfs$, a distinguishing feature of the solution decomposed approach is that both sets of parameters $\bfs_1$ and $\bfs_2$ are estimated from the data $\bfd$, even if $\bfs$ is desired in the end.

\paragraph{Contributions and overview}  In this paper we propose a new computational framework for solving inverse problems in which the solution is assumed to be the sum of two components: a ``smooth'' background and a sparse term that represents anomalies. Following the Bayesian approach, we derive a posterior distribution for the unknown solution components and focus on efficiently computing the MAP estimate. Our approach has three main components:
\begin{enumerate}
    \item  Use of a majorization-minimization (MM) scheme to solve the optimization problem for the MAP estimate as a sequence of iteratively reweighted least-squares problems that carefully reweights only the sparse term, 
    \item A novel Flexible, Generalized Golub-Kahan (FGGK) iterative method for efficiently generating a single basis for approximately solving the reweighted least-squares subproblems, and 
    \item Robust methods for automatically selecting the regularization parameters within the projected solution space at each iteration. 
\end{enumerate}
The main novelty of this paper is that in contrast to inner-outer methods or alternating approaches that solve a sequence of least-squares/optimization problems from ``scratch,'' our approach successively builds a single basis with which we seek approximate solutions for the components. To accomplish this task, we develop the FGGK approach, a new Krylov subspace solver, by paying close attention to the computational cost; like the generalized Golub-Kahan (genGK) process \cite{chung2017generalized}, it avoids forming inverses of covariance matrices and uses the same number of matrix-vector products with the forward operator and covariance matrices, but has the ability to incorporate information about both the solution components. By automatically selecting the regularization parameter at each iteration, we avoid the need to solve repeated optimization problems to estimate the parameters. 
 
 We demonstrate the performance of our approach on a series of large-scale test problems that represent anomaly detection in dynamic tomography and atmospheric inverse modeling, where for the last example that contains over a million unknowns, our method is capable of achieving satisfactory results in fewer than $50$ iterations.

An outline of the paper is as follows.  In \cref{sec:background} we describe a Bayesian approach for \cref{eq:sdproblem} and provide an overview of related works that incorporate multiple stochastic components for inverse problems.
The proposed hybrid projection method for solution decomposition is described in \cref{sec:hybridmethod}, along with a description of some algorithmic considerations as well as methods to choose regularization parameters. 
Numerical results are provided in \cref{sec:numerics}, and conclusions are provided in \cref{sec:conclusions}.

\section{Bayesian inverse problems with solution decomposition}
\label{sec:background}

In this section, we describe a Bayesian approach for solving \cref{eq:sdproblem}.
We assume that the errors $\bfdelta$ and the unknowns $\bfs$ are mutually independent.
For anomaly detection and more generally for problems where the solution can be decomposed, we consider the case where $\bfs_1$ follows a multivariate Gaussian distribution, and the components of $\bfs_2$ are independent and follow the univariate Laplace distribution~\cite[Chapter 4.3]{bardsley2018computational}. The Laplace prior enforces sparsity in $\bfs_2$ (e.g., the solution itself is sparse and contains many zeros or a representation in some frequency domain is sparse) \cite{kaipio2006statistical}.

In the geostatistical framework~\cite{matheron_1973,kitanidis1983geostatistical,kitanidis1986parameter,kitanidis1995quasi}, we model the unknown function $s(\bfzeta)$, where $\bfzeta \in \mathbb{R}^d$ represents the coordinates in space, as a realization of a random field. We express this realization as a sum of two terms:
\[ s(\bfzeta) = s_1(\bfzeta) + s_2(\bfzeta) \qquad s_2(\bfzeta) := \sum_{k=1}^p \beta_k \psi_k(\bfzeta) , \]
where $s_1$ is the realization of a random field that captures the smooth features, $\psi_k$ are deterministic basis functions, and $\beta_k$ are coefficients to be determined. As is the prevalent approach, we take the random field $s_1$ to be Gaussian, which is characterized by a mean function $\mu_1(\bfzeta)$ and a covariance function $\lambda^{-2}\kappa(\bfzeta,\bfzeta')$, where $\lambda$ is a parameter that controls the precision and is a hyperparameter that must be determined.  For short, we write $s_1 \sim \text{GP}(\mu_1(\bfzeta),\lambda^{-2}\kappa(\bfzeta,\bfzeta'))$, where GP denotes Gaussian process. We consider a set of grid points $\{\bfzeta_j \}_{j=1}^n$ on which we represent the unknown random field. Define the vector $\bfs_1 = \begin{bmatrix} s_1(\bfzeta_1) & \dots & s_1(\bfzeta_n) \end{bmatrix}^\top$, then it follows:
\[ \bfs_1 \sim \calN(\bfmu_1, \lambda^{-2} \bfQ),\]
where $\bfmu_1 = \begin{bmatrix} \mu_1(\bfzeta_1) & \dots & \mu_1(\bfzeta_n) \end{bmatrix}^\top$ and $\bfQ_{ij} = \kappa(\bfzeta_i,\bfzeta_j)$ for $1\leq i, j \leq n$.

In this paper, to model the anomalies, we use point sources  whose locations coincide with the grid points; this can be accomplished by taking the basis functions $\psi_k 
(\bfzeta) = \delta(\bfzeta- \bfzeta_k) $ as Dirac delta functions and the number of basis functions $p = n$. In practice, we represent the delta function by a Gaussian with a scale parameter smaller than the mesh width.  Since anomalies are localized, we enforce sparsity in the coefficients $\beta_k$ as follows: we assume that the coefficients $\beta_k$ are independent of one another and the random field $s_1$, and impose the univariate Laplace distribution with mean $[\bfmu_2]_j$ and scale parameter $2\alpha^{-2}$, that is, 
\[ \beta_j \sim \mathcal{L}([\bfmu_2]_j, 2\alpha^{-2}) \qquad 1 \leq j \leq n.\]
Similar to $\bfs_1$, we define the vector $\bfs_2$ with components $[\bfs_2]_j = \beta_j $ for $1 \leq j \leq n$.

In summary, we have the prior model
\begin{equation}
    \label{eq:assumptions}
    \bfs_1 \sim \calN(\bfmu_1, \lambda^{-2} \bfQ) \quad \mbox{and} \quad [\bfs_2]_j \sim \calL ([\bfmu_2]_j, 2 \alpha^{-2}), \quad 1 \leq j \leq n,
\end{equation} where $\bfmu_1, \bfmu_2 \in \bbR^n$, $\bfQ$ is SPD, and $\lambda \neq 0, \alpha \neq 0$ are scaling parameters.
For problems of interest, computing the inverse and square root of $\bfR$ is inexpensive (e.g., $\bfR$ is often a diagonal matrix), but explicit computation of $\bfQ$ (or its inverse or square root) may not be possible.  However, we assume that matrix-vector multiplications (matvecs) involving $\bfA$, $\bfA\t$, and $\bfQ$ can be done efficiently (e.g., in $\calO(n \log n)$ operations rather than $\calO(n^2)$ operations for an $n \times n$ matrix); for details, see~\cite{ambikasaran2012large}. This framework can be extended to spatiotemporal models,  but we do not provide the details here.

Under assumptions \cref{eq:sdproblem,eq:assumptions} and using Bayes' theorem, the posterior probability density function is given by
\begin{equation} \label{eq:posterior}
\begin{aligned}
    \pi_{\rm post} (\bfs_1,\bfs_2 \mid \bfd) &= \frac{\pi(\bfd \mid \bfs_1,\bfs_2) \pi(\bfs_1)\pi(\bfs_2)}{\pi(\bfd)} \\
    & \propto \exp\left( -\frac{1}{2}\norm[{\bfR^{-1}}]{\bfA \bfs - \bfd}^2 - \frac{\lambda^2}{2} \norm[\bfQ^{-1}]{\bfs_1-\bfmu_1}^2  - \frac{\alpha^2}{2} \norm[1]{\bfs_2 - \bfmu_2}\right),
\end{aligned}
\end{equation}
where $\|\cdot\|_1$ denotes the 1-norm of a vector, $\norm[\bfM]{\bfx}^2 = \bfx\t \bfM \bfx$ for any SPD matrix $\bfM$, and $\propto$ means ``proportional to.'' In the Bayesian framework, the solution is the posterior distribution.  Notice that the posterior is not Gaussian. 

In this manuscript, we 
describe new hybrid projection methods to efficiently approximate the maximum a posteriori (MAP) estimate, which corresponds to the mode of the posterior distribution and is the solution to the following optimization problem,
\begin{equation}
\label{eq:map}
	\min_{\bfs_1 \in \R^n, \bfs_2\in \R^n} \norm[{\bfR^{-1}}]{\bfA (\bfs_1+\bfs_2) - \bfd}^2 + \lambda^2 \norm[\bfQ^{-1}]{\bfs_1-\bfmu_1}^2 + \alpha^2 \norm[1]{\bfs_2 - \bfmu_2}.
\end{equation}

\paragraph{Computational challenges}Optimization problems such as \eqref{eq:map} can be computationally challenging to solve, and there are three main computational concerns.  First, an accurate reconstruction of $\bfs$ will rely heavily on being able to obtain good estimates of regularization parameters $\lambda$ and $\alpha$, which can be very difficult to estimate prior to solution computation.  Second, for many problems with nonstandard priors (e.g., priors defined on nonstructured grids), explicit computation of $\bfQ$ (or its inverse or square root) may not be possible.  Generalized hybrid iterative methods which are based on the genGK bidiagonalization can be used to solve problems of the form
\begin{equation}
\label{eq:genmap}
	\min_{\bfs_1 \in \R^n} \norm[{\bfR^{-1}}]{\bfA \bfs_1 - \bfd}^2 + \lambda^2 \norm[\bfQ^{-1}]{\bfs_1-\bfmu_1}^2
\end{equation}
and are described in~\cite{chung2017generalized}.  Third, it is well known that solving the $\ell_1$ regularized problem can be computationally difficult, due to nondifferentiability at the origin as well as the need to use expensive nonlinear or iteratively reweighted optimization schemes.
These inner-outer approaches can get very costly \cite{fornasier2016cg, wohlberg2008lp}, which has led to the development of accelerated alternative methods such as split Bregman methods \cite{xiong2019convex} and iterative shrinkage threshholding algorithms \cite{beck2009fast}, where an iterative two-step process is used.  However, these methods require a priori selection of various parameters (including the regularization and shrinkage parameter), which can be cumbersome \cite{matsui2011variable}.  Flexible Krylov methods have been proposed \cite{gazzola2014generalized,chung2019flexible} as a means to both avoid inner-outer schemes and allow automatic regularization parameter selection.  Such methods can be used to solve problems of the form,
\begin{equation}
\label{eq:flexmap}
	\min_{\bfs_2 \in \R^n} \norm[2]{\bfA \bfs_2 - \bfd}^2  + \alpha^2 \norm[1]{\bfs_2-\bfmu_2}
\end{equation}
and are described in~\cite{chung2019flexible}. 

A natural approach to solve \cref{eq:map} would utilize an alternating optimization scheme (see \cref{alg:altopt}), but this approach may be very slow and requires unrealistic initialization vectors.  Instead, we describe in \cref{sec:hybridmethod} an iterative FGGK approach that combines the flexible and generalized Golub-Kahan projection methods to efficiently generate a basis for solving inverse problems with mixed prior models.

\begin{algorithm}[!ht]
\begin{algorithmic}[1]
\STATE Initialize $\bfs_2^{(0)}, k = 1$
\WHILE {not converged}
\STATE Solve $\bfs_1^{(k)} = \argmin_{\bfs_1} \norm[{\bfR^{-1}}]{\bfA \bfs_1 - (\bfd-\bfA \bfs_2^{(k-1)})}^2 + \lambda^2 \norm[\bfQ^{-1}]{\bfs_1-\bfmu_1}^2$
\STATE Solve
$\bfs_2^{(k)} = \arg\min_{\bfs_2} \norm[2]{\bfA \bfs_2 - (\bfd-\bfA \bfs_1^{(k)})}^2  + \alpha^2 \norm[1]{\bfs_2-\bfmu_2}$
\STATE Set $k = k+1$
\ENDWHILE
\end{algorithmic}
\caption{Solving \cref{eq:map} using alternating optimization}
\label{alg:altopt}
\end{algorithm}

\paragraph{Related approaches and ideas}
The idea to incorporate multiple stochastic components for inverse problems is not necessarily new, nor is it restricted to atmospheric imaging. For example, different texture models have been investigated for improved breast cancer imaging \cite{li2016novel}.  However, previous methods to handle multiple stochastic components are quite costly and often rely on simplifying assumptions.  For example, in 
\cite{yadav2016statistical}, the authors disaggregate the unknown fluxes to account for the biospheric and fossil fuel components separately, but simple Gaussian priors were used (a stationary, separable exponential model for the biospheric fluxes and a diagonal covariance matrix for the fossil fuel fluxes). In the special case where both priors are Gaussian, i.e., $\bfs_1 \sim \calN(\bfmu_1, \bfQ_1)$ and $\bfs_2 \sim \calN(\bfmu_2, \bfQ_2)$ with $\bfmu_1, \bfmu_2 \in \bbR^n$ and SPD matrices $\bfQ_1, \bfQ_2$, it can be shown that $\bfs \sim \calN(\bfmu_1+\bfmu_2, \bfQ_1+\bfQ_2)$, and methods for mixed Gaussian priors can be used \cite{cho2020hybrid}. 
Furthermore, if sparsity or a sparsity decomposition is desired, methods based on robust PCA have been developed, e.g. in dynamic magnetic resonance imaging, to separate the solution (reshaped into a matrix) into a low-rank plus a sparse matrix \cite{tremoulheac2014dynamic}, but such methods are too restrictive for the problems of interest.

We remark that although optimization problem \cref{eq:map} has a similar flavor to elastic net regularization \cite{zou2005regularization,jin2009elastic,chen2016elastic} and other $\ell_1-\ell_2$ problems \cite{yin2015minimization}, our approach is fundamentally different. 
First, these approaches linearly combine $\ell_1$ and $\ell_2$ regularization for $\bfs$ and do not split the solution into two stochastic components.  That is, they assume that the entire solution vector is included in both regularization terms, e.g., for elastic net,
\begin{equation}
    \label{eq:elasticnet}
    \min_{\bfs \in \R^n} \norm[2]{\bfA \bfs - \bfd}^2 + \lambda^2 \norm[2]{\bfs}^2 + \alpha^2 \norm[1]{\bfs}.
\end{equation}
Second, typically iterative numerical methods of active set type are employed to solve elastic net regularized problems, but the inclusion of nontrivial regularizers makes these approaches infeasible.

\section{Hybrid projection methods for solution decomposition}
\label{sec:hybridmethod}
In this section, we describe an efficient computational method to approximate the MAP estimate given in \eqref{eq:map}. We develop a \textit{combined} hybrid projection method that builds on the generalized and flexible Golub-Kahan processes and inherits many of the main computational benefits from previously developed hybrid approaches.

We begin in \cref{sub:mmapproach} by describing an MM approach to handle the $1$-norm regularization term in \eqref{eq:map}. Although this inner-outer optimization approach is computationally infeasible, it motivates the use of flexible preconditioning techniques that, extended and combined with genGK methods, are described in \cref{ss_hybridmethod}.  In particular, we describe a Flexible Generalized Golub-Kahan (FGGK) projection method, with some discussion on algorithmic considerations (e.g., breakdown) and computational considerations, and then we describe a hybrid projection approach based on the FGGK projection. We pay special attention to methods to choose regularization parameters for the projected problem in \cref{ss_regpar}.

\subsection{Majorization-Minimization (MM) approach} 
\label{sub:mmapproach}
Various methods have been developed for approximating the solution of $\ell_1$-regularized problem \cref{eq:flexmap}, ranging from iterative shrinkage algorithms to iterative reweighted norms \cite{beck2009fast,wohlberg2008lp,IRNekki}. In this subsection, we provide an overview of the MM approach for approximating the solution of \cref{eq:map}, which requires solving a sequence of optimization problems.  For this discussion, we assume that $\lambda$ and $\alpha$ are fixed.

We begin with the following change of variables, 
\begin{equation}\label{eqn:changeofvar}
    \bfs_1 = \bfmu_1 + \bfQ \bfx, \quad \bfs_2 = \bfmu_2 + \bfxi, \quad \mbox{and} \quad \bfc = \bfd - \bfA \bfmu_1 - \bfA \bfmu_2,
\end{equation} 
and get optimization problem
\begin{equation}\label{eq:mix3}
	\min_{\bfx \in \bbR^n,\, \bfxi \in \bbR^n} f(\bfx,\bfxi) =  \norm[\bfR^{-1}]{\bfA \bfQ \bfx + \bfA \bfxi - \bfc}^2 + \lambda^2 \norm[\bfQ]{\bfx}^2 + \alpha^2 \norm[1]{\bfxi}.
\end{equation}
Notice that we have removed all instances of $\bfQ^{-1}$.  Next to handle the $\ell_1$ term, we use the MM approach to convert optimization problem \cref{eq:mix3} into a sequence of reweighted least-squares problem. 
For some $\epsilon > 0$, we consider $|t| \approx \phi_\epsilon(t) = \sqrt{t^2 + \epsilon}$ and approximate $\|\bfxi\|_1 \approx \sum_{j=1}^n\phi_\epsilon(\xi_j)$. The corresponding objective function is
\begin{align}
\label{eq:objfun}
    f_\epsilon(\bfx,\bfxi) = \norm[\bfR^{-1}]{\bfA \bfQ \bfx + \bfA \bfxi - \bfc}^2 + \lambda^2 \norm[\bfQ]{\bfx}^2 + \alpha^2\sum_{j=1}^n\phi_\epsilon( \xi_j).
\end{align}
From Equation 1.5 in \cite{lange2016mm}, we have the majorization relationship
\[ \phi_\epsilon(t) = \sqrt{t^2 + \epsilon} \leq \sqrt{t_k^2 + \epsilon } + \frac{1}{2\sqrt{t_k^2 + \epsilon}}(t^2 - t_k^2) =: \psi_\epsilon(t \mid t_k).\]
Given the current iterate $(\bfx^{(k)}, \bfxi^{(k)}),$ we can define the surrogate function
\[ g_\epsilon(\bfx,\bfxi \mid \bfx^{(k)},\bfxi^{(k)}) = \norm[\bfR^{-1}]{\bfA \bfQ \bfx + \bfA \bfxi - \bfc}^2 + \lambda^2 \norm[\bfQ]{\bfx}^2 + \alpha^2\sum_{j=1}^n \psi_\epsilon(\xi_j \mid \xi^{(k)}_j). \]
It is easily verified that $f_\epsilon(\bfx^{(k)},\bfxi^{(k)}) = g_\epsilon(\bfx^{(k)},\bfxi^{(k)} \mid \bfx^{(k)},\bfxi^{(k)})$ and 
\[ f_\epsilon(\bfx,\bfxi) \leq  g_\epsilon(\bfx^{(k)},\bfxi^{(k)} \mid \bfx^{(k)},\bfxi^{(k)}) \qquad \forall  \bfx, \bfxi  \in \bbR^n. \]
These two conditions mean that the surrogate function $g_\epsilon(\bfx,\bfxi \mid \bfx^{(k)})$ matches the objective function $f_\epsilon(\bfx,\bfxi)$ at the current iterate and majorizes the surrogate function for every point, respectively. This means that as long as we choose the next iterate $(\bfx^{(k+1)},\bfxi^{(k+1)})$ such that the surrogate $g_\epsilon$ decreases, then we ensure that this decreases the objective $f_\epsilon$ since
\[ \begin{aligned} f_\epsilon(\bfx^{(k+1)},\bfxi^{(k+1)}) \leq & \>  g_\epsilon(\bfx^{(k+1)},\bfxi^{(k+1)}\mid \bfx^{(k)},\bfxi^{(k)}) \\ 
 \leq & \> g_\epsilon(\bfx^{(k)},\bfxi^{(k)}\mid\bfx^{(k)},\bfxi^{(k)}) = f_\epsilon(\bfx^{(k)},\bfxi^{(k)}). \end{aligned}\]
 The first inequality and the final equality are due to the majorization properties of the surrogate function $g_\epsilon$. The second inequality is satisfied by choosing the next iterates in a manner to ensure that the surrogate is decreased; it is important to note that it is not necessary to minimize the surrogate at each iteration. 
 
 Thus, the MM algorithm for solving~\eqref{eq:mix3} is as follows: Given initial guesses $(\bfx^{(0)},\bfxi^{(0)})$, solve the following sequence of reweighted least-squares problems
 \begin{align}(\bfx^{(k+1)},\bfxi^{(k+1)}) & = \argmin_{\bfx\in \R^n ,\bfxi\in \R^n} g_\epsilon(\bfx,\bfxi \mid \bfx^{(k)},\bfxi^{(k)}) \nonumber \\ 
 & = \argmin_{\bfx \in \R^n, \bfxi \in \R^n} \norm[\bfR^{-1}]{\bfA \bfQ \bfx + \bfA \bfxi - \bfc}^2 + \lambda^2 \norm[\bfQ]{\bfx}^2 + \alpha^2 \norm[2]{\bfD(\bfxi^{(k)})\bfxi}^2, 
 \label{eq:MM_min}
 \end{align}
where terms from $g_\epsilon$ that do not depend on $\bfx$ and $\bfxi$ have been dropped
 and diagonal matrix $\bfD(\bfxi) = \text{diag}\left( \left[2\sqrt{\xi_i^2 + \epsilon}\right]^{-1/2}\right)_{i=1}^n.$ 
     \label{eq:MMhat_min}
To get the solution after removing the change of variables, we get $\bfs^{(k+1)} = \bfmu_1 + \bfQ\bfx^{(k+1)} + \bfmu_2 +  \bfxi^{(k+1)}$.
  
 The convergence of the MM scheme has been established, see e.g., \cite{huang2017some}.  However, notice that minimizing the surrogate requires solving a large optimization problem \cref{eq:MMhat_min} with $2n$ unknowns at each iteration. For small problems, one could solve the corresponding normal equations. For larger problems, an iterative method could be used to solve the reweighted least-squares problems, but this often leads to expensive inner-outer solves \cite{wohlberg2008lp}.  Instead, we describe in the next section an approach that avoids inner-outer schemes by exploiting flexible preconditioning techniques, following recent works, e.g., \cite{gazzola2014generalized,chung2019flexible}.  That is, 
 we solve optimization problem~\eqref{eq:MMhat_min} approximately at each step.

\subsection{Flexible Generalized Golub-Kahan (FGGK) iterative method}
\label{ss_hybridmethod} 
In this section, we describe iterative projection methods that can be used for approximating the solution for inverse problems with solution decomposition (e.g., for anomaly detection).
We exploit aspects of both the flexible and generalized Golub-Kahan projection methods and develop a solution decomposition hybrid projection approach, henceforth dubbed \texttt{sdHybr}, to approximate the MAP estimate \cref{eq:map}.  Similar to all hybrid projection methods that combine iterative projection methods with variational regularization techniques, there are two main components.  First, we generate a basis for the solution (which in this case includes two sets of solution vectors) by exploiting a flexible preconditioning framework integrated with a genGK bidiagonalization.  Second, we compute an approximate solution to the inverse problem by solving an optimization problem in the projected subspace where regularization parameters can be estimated automatically.
  
\paragraph{FGGK process} We consider solving problems of the form \cref{eq:MMhat_min} by incorporating a changing diagonal matrix, which rescales the norms, to generate a basis for the solution.  Suppose we are given $\bfA$, $\bfQ,$ $\bfR,$ $\bfc$, and a sequence of invertible matrices  $\{\bfD_j\}_{j=1}^k.$ We initialize the iterations with  $m_{1,1} = \|\bfc\|_{\bfR^{-1}}$ and $\bfu_1 = \bfc/m_{1,1}$; furthermore, take $\bfv_1 = \bfA^\top\bfR^{-1}\bfu_1$ and $t_{1,1} = \| \bfv_1\|_{\bfQ} $. The FGGK iterative process generates vectors $\bfz_k,$ $\bfv_k$, and $\bfu_{k+1}$ such that at the $k$th iteration, we have 
\begin{eqnarray}
\label{iter1}
    m_{k+1,k} \bfu_{k+1} = & \>  \bfA \bfQ \bfv_k  + \bfA \bfD_k^{-1}\bfv_k - \sum_{j=1}^k m_{j,k}\bfu_j \\
\label{iter2}
    t_{k+1,k} \bfv_{k+1} = & \>  \bfA^\top \bfR^{-1}\bfu_{k+1} - \sum_{j=1}^k t_{j,k}\bfv_j. 
\end{eqnarray}
In the first step, we expand the basis $\bfu_{k+1}$ by including the vectors $\bfA \bfQ \bfv_k$ and $\bfA \bfD_k^{-1}\bfv_k$, and we orthogonalize against the previous basis vectors $\bfu_1,\dots,\bfu_k$ using Gram-Schmidt with the inner product $\langle \cdot,\cdot\rangle_{\bfQ}$. Similarly, we expand the basis vectors $\bfv_{k+1}$ with the vector $\bfA^\top \bfR^{-1}\bfu_{k+1}$ and orthogonalize against the previous vectors $\bfv_1,\dots,\bfv_k$ using the inner product $\langle \cdot,\cdot\rangle_{\bfR^{-1}}$. Finally we ensure that both $\bfu_{j},\bfv_{j}$ are normalized so that $\|\bfu_{j}\|_{\bfR^{-1}} =\|\bfv_{j}\|_{\bfQ} = 1$, and let $\bfw_j = \bfD_j^{-1} \bfv_j$ for $1 \leq j \leq k+1$. 

For notational convenience, consider the augmented matrices  
\begin{equation}\label{eqn:qhat}
   \widehat \bfA = \begin{bmatrix} \bfA & \bfA \end{bmatrix} \in \bbR^{m \times 2n} \qquad \mbox{and} \qquad \widehat \bfQ = \begin{bmatrix} \bfQ & \\ & \bfI \end{bmatrix} \in \bbR^{(2n) \times (2n)}.
\end{equation} The equations~\eqref{iter1} and~\eqref{iter2} can be summarized in matrix form as 
\begin{equation}
\label{eq:genflexrel}
	\widehat \bfA \widehat \bfQ \bfZ_k = \bfU_{k+1} \bfM_k \quad \mbox{and} \quad \bfA\t \bfR^{-1} \bfU_{k+1} = \bfV_{k+1} \bfT_{k+1},
\end{equation}
where the search basis $\bfZ_k$ takes the form
\begin{equation} 
\label{eq:Zk}\bfZ_k = \begin{bmatrix} \bfz_1 & \dots & \bfz_k \end{bmatrix} = \begin{bmatrix}  \bfv_1 & \dots &  \bfv_k \\ \bfw_1 & \dots & \bfw_k \end{bmatrix} = \begin{bmatrix}  \bfV_k \\ \bfW_k \end{bmatrix}\in \bbR^{2n \times k}.
\end{equation}
Furthermore, we have two matrices $\bfM_k \in \bbR^{(k+1)\times k}$ and $\bfT_{k+1} \in \bbR^{(k+1)\times (k+1)}$ that are upper Hessenberg and upper triangular respectively.
The basis vectors are collected in the matrices $\bfU_{k+1} =
\begin{bmatrix}
 \bfu_1 & \dots & \bfu_{k+1}
\end{bmatrix} \in \bbR^{m\times (k+1)}$ and $\bfV_{k+1} =\begin{bmatrix}
 \bfv_1 & \dots & \bfv_{k+1}
\end{bmatrix}\in \bbR^{n\times (k+1)}$ which satisfy the orthogonality conditions (in exact arithmetic)
\begin{equation}
    \label{eq:QV-orth}
    \bfU_{k+1}\t \bfR^{-1}\bfU_{k+1} = \bfI_{k+1} \qquad \bfV_{k+1}\t \bfQ \bfV_{k+1} =\bfI_{k+1}.
\end{equation}

\paragraph{Solving the LS problem} Thus far, we have described the FGGK process as an iterative method to generate a basis for the solution.  Next, we seek approximate solutions to the least-squares problem \cref{eq:MMhat_min} by using the FGGK relations above to obtain a sequence of projected LS problems. For clarity of presentation, we assume that the parameters $\lambda$ and $\alpha$ are fixed.   

To determine the optimal coefficients $\bff$ in~\eqref{eqn:sdhybrmin}, we plug the FGGK relations~\eqref{eq:genflexrel} into the objective function to get 
\begin{equation}
\label{eq:projectedFGGK}
	\bff_k = \argmin_{\bff \in \bbR^k} \norm[2]{\bfM_k \bff - m_{1,1}\bfe_1}^2 +  \lambda^2 \norm[2]{\bff}^2 +\alpha^2 \norm[2]{\bfW_k \bff}^2.
	\end{equation}	
which is equivalent to
\begin{equation}\label{eqn:sdhybrmin} \min_{\bfx =\bfV_k\bff, \bfy = \bfW_k\bff} \| \bfA \bfQ \bfx + \bfA\bfy - \bfc \|_{\bfR^{-1}}^2 + \lambda^2 \|\bfx \|_{\bfQ}^2 + \alpha^2\|\bfy\|_2^2. \end{equation}

Our ultimate goal is to obtain the decomposed solution $\bfs= \bfs_1  + \bfs_2$. Based on the discussion above and using the change of variables~\eqref{eqn:changeofvar}, we get  
\[ \bfs_1^{(k)} := \bfmu_1 + \bfQ \bfV_k \bff_k  \qquad \bfs_2^{(k)} := \bfmu_2 + \bfW_k\bff_k ,  \]
where $\bff_k \in \bbR^k$ are the coefficients obtained by solving~\eqref{eq:projectedFGGK}.
In summary, an approximation to the MAP estimate at the $k$th iteration of the \texttt{sdHybr} method would be $\bfs^{(k)} = \bfmu_1 + \bfQ \bfV_k \bff_k + \bfmu_2 + \bfW_k \bff_k$ where 
$\bff_k$ solves the optimization problem~\eqref{eq:projectedFGGK}.

\paragraph{Efficient QR update} Similar to what is done in \cite{chung2019flexible}, we can take the thin QR factorization of $\bfW_k = \bfQ_{W,k}\bfR_{W,k}$, then~\eqref{eq:projectedFGGK} becomes
\begin{equation}
\label{eq:projected}
    \bff_k = \argmin_{\bff \in \bbR^k} \norm[2]{\bfM_k \bff - m_{1,1}\bfe_1}^2 +  \lambda^2 \norm[2]{\bff}^2 +\alpha^2 \norm[2]{{\bfR}_{W,k} \bff}^2. 
\end{equation}
An efficiently update of the QR factorization of $\bfW_k$ is as follows.  Suppose we have have computed the thin-QR factorization 
\begin{equation}
    \bfW_k = \bfQ_{W,k}\bfR_{W,k}
\end{equation}
 where $\bfQ_{W,k}^{\top}\bfQ_{W,k} = \bfI_{k}$ and $\bfR_{W,k}$ is an upper triangular matrix. Then we can efficiently update the QR factorization of $\bfW_k$ in $\mathcal{O}(nk)$ flops using Gram-Schmidt as 
\begin{equation}
\label{qr1}
\begin{aligned}
    \bfW_{k+1} &= \begin{bmatrix}
     \bfW_k & \bfw_{k+1}
    \end{bmatrix} \\
    & = \begin{bmatrix}
     \bfQ_{W,k}\bfR_{W,k} &  \bfw_{k+1}
    \end{bmatrix} \\
    &= \underbrace{\begin{bmatrix}
     \bfQ_{W,k} & (\bfI - \bfQ_{W,k}\bfQ_{W,k}^{\top}) \bfw_{k+1}/\beta_{k+1}
    \end{bmatrix}}_{\bfQ_{W,k+1}}\underbrace{\begin{bmatrix}
     \bfR_{W,k} & \bfQ_{W,k}^{\top} \bfw_{k+1} \\
     \bf0 & \beta_{k+1}
    \end{bmatrix}}_{\bfR_{W,k+1}}, 
\end{aligned}
\end{equation}
where $\beta_{k+1} = \norm[2]{(\bfI - \bfQ_{W,k}\bfQ_{W,k}^{\top}) \bfw_{k+1}}$. 
    For additional numerical stability, one can use another round of Gram-Schmidt or instead use Householder QR updates.
A summary of the \texttt{sdHybr} method
is provided in \cref{alg:hybridmethod}.  Parameter selection methods to select $\alpha$ and $\lambda$ automatically will be described in \cref{ss_regpar}.  We remark that various existing stopping criteria can be used and will be discussed shortly.

\paragraph{Computational Cost} Each iteration of the \texttt{sdHybr} method requires one matrix-vector products with $\bfA$ and its adjoint (suppose we denote the cost of one matrix-vector of this operation by $T_{\bfA}$), two matrix-vector products with $\bfQ$ (similarly, denoted $T_{\bfQ}$), one matrix-vector product  with $\bfR^{-1}$ (denoted $\bfT_{\bfR^{-1}}$), one matrix-vector product with $\bfD_k^{-1}$ (denoted $T_{\bfD_k^{-1}}$), the inversion of diagonal matrix $\bfD_k$ that is $\mathcal{O}(n)$ floating point operations (flops),
and additional $\mathcal{O}(k(m+n))$ flops for the summation calculation in \eqref{iter1} and \eqref{iter2}. To compute the solution of the projected problem \eqref{eq:projected}, the cost is $\mathcal{O}(k^3)$ flops, since $\bfM_k$ is upper Hessenberg matrix.  And the cost of forming $\bfx$ and $\bfy$ to get $\bfs$ is $\mathcal{O}(k^2(m+n))$ flops. 
Thus, the overall cost of the algorithm is 
\begin{equation}
    T_{\text{sdHybr}} = 2kT_{\bfA} + 2kT_{\bfQ} + kT_{\bfR^{-1}} + \mathcal{O}(k^2(m+n))  + \mathcal{O}(k^4)   
    \; \text{flops}.
\end{equation}
\begin{algorithm}[!ht]
\caption{Solution decomposition hybrid method (\texttt{sdHybr})}
\label{alg:hybridmethod}
\begin{algorithmic}[1]
	\REQUIRE{Matrix $\bfA \in \mathbb{R}^{m \times n}$, positive definite matrices $\bfQ \in \mathbb{R}^{n \times n}$ and $\bfR \in \mathbb{R}^{m \times m}$, vector $\bfc \in \mathbb{R}^{m}, \bfmu_1, \bfmu_2 \in \mathbb{R}^{n} $.  Invertible matrix $\bfD_1 = \bfI_n \in \mathbb{R}^{n \times n}$. } 
\STATE Initialize $\bfu_1 = \bfc/m_{1,1}$, where $m_{1,1} = \norm[\bfR^{-1}]{\bfc}$ and $\bfv_1 = \bfzero, k = 1$. 
\WHILE{stopping criteria are not satisfied}
\STATE $\bfw = \bfA\t \bfR^{-1} \bfu_k$, $t_{j,k} = \bfw\t  \bfQ \bfv_j$ for $j=1, \ldots, k-1$
\STATE $\bfh = \bfh - \sum_{j=1}^{k-1} t_{j,k}  \bfv_j$, $t_{k,k} = \norm[\bfQ]{\bfh}$, $\bfv_k = \bfh/t_{k,k}$
\STATE $\bfz_k = \begin{bmatrix} \bfv_k \\ \bfw_k \end{bmatrix} $, $\bfV_k = \begin{bmatrix}
   \bfv_1 & \dots & \bfv_k 
  \end{bmatrix}$, $\bfW_k = \begin{bmatrix}
   \bfw_1 & \dots & \bfw_k 
  \end{bmatrix}$, where $\bfw_k = \bfD_k^{-1} \bfv_k$.
\STATE 
$\bfh = \bfA(\bfQ\bfv_k + \bfw_k)$, $m_{j,k} = \bfh\t \bfR^{-1}\bfu_j$ for $j=1, \ldots, k$
\STATE $\bfh = \bfh - \sum_{j=1}^k m_{j,k} \bfu_j$, $m_{k+1,k} = \norm[\bfR^{-1}]{\bfh}$, $\bfu_{k+1} = \bfh/m_{k+1,k}$
\STATE Update QR factorization using \eqref{qr1} to obtain $\bfW_{k+1} = \bfQ_{k+1}\bfR_{k+1}$
\STATE Solve \cref{eq:projected} to get $\bff_k(\lambda_k,\alpha_k)$ with selected regularization parameters $\lambda_k,\alpha_k$.
\STATE $\bfs_1^{(k)} = \bfmu_1 + \bfQ\bfV_k\bff_k$, $\bfs_2^{(k)} = \bfmu_2 + \bfW_k\bff_k$.
  \STATE $\bfD_{k+1} =  \bfD(\bfW_k\bff_k)$
 \STATE $k = k+1$
\ENDWHILE
	\RETURN{Approximations $\bfs_1^{(k)}$ and $\bfs_2^{(k)}$ that define the sum $\bfs^{(k)} =\bfs_1^{(k)} + \bfs_2^{(k)}$.}
\end{algorithmic}
\end{algorithm}

Notice that compared to the MM method, the projected problem~\eqref{eq:projectedFGGK} for \texttt{sdHybr} is cheaper to solve than~\eqref{eq:MM_min} at each MM iteration, since it is an optimization problem over a smaller dimensional space.

\subsubsection*{Alternative approaches} We remark that various projection methods that combine flexible and generalized GK methods could be considered besides the \texttt{sdHybr} approach; however, for stability and for proper selection of regularization parameters, we found this projection approach (\cref{alg:hybridmethod}) to be the most computationally appealing. For example, a na\"ive first approach would be to use a genGK approach to handle the $\bfQ$-norm regularizer and use the flexible Golub-Kahan approach to handle the $1$-norm separately in \cref{eq:mix3}. This would generate two solution subspaces that each contain orthonormal columns but are not orthogonal to each other, and the number of unknowns for the projected problem would be $2k$. Although an efficient QR update could be used, there are potential issues with breakdown (e.g., due to linear dependence of subspace vectors).  Furthermore, we found that this approach can be sensitive to initializations. The FGGK process described above avoids this by working with stacked solution vectors and a projected problem of order $k$.  

Another natural approach to combine flexible and generalized GK methods would be to reformulate the problem such that the genGK vectors include the flexible preconditioner. This approach is described in \cref{subsec:app1} and provides a nice alternative, but the main caveat is that selecting regularization parameters in a hybrid framework becomes more challenging.  In particular, a hybrid framework based on the projection method described in \cref{subsec:app1} requires formulating one regularization parameter as a fixed scalar multiple of the other or utilizing more expensive optimization procedures.  This is due to the inability to simplify the norm for the regularizer when different regularization parameters. In contrast, the FGGK procedure leads to a projected problem \cref{eq:projectedFGGK} for which existing parameter selection approaches can be naturally applied, as described in the next section.

\subsubsection*{Remarks on solution decomposition method}
The FGGK process was designed to solve the sequence of optimization problems~\eqref{eq:MM_min} involving a diagonal weighting matrix that changes at each iteration. We make a few comments when the weighting matrix is fixed at each iteration.
Recall that the  Krylov subspace associated with the matrix $\bfE$ and vector $\bfg$ is defined as 
\begin{equation*}
    \mathcal{K}_k(\bfE, \bfg) \equiv \text{Span}\{\bfg, \bfE\bfg, \dots, \bfE^{k-1}\bfg\}.
\end{equation*}
Assume that  $\bfD_j$ with $j = 1, \dots, k$ are fixed, that is $\bfD_1 = \dots = \bfD_k = \widehat{\bfD}$, and furthermore, assume $\widehat{\bfD}$ is invertible.
For the FGGK process, based on the relations in \eqref{iter1} and \eqref{iter2}, it can be shown that the columns of $\bfU_k$ $\bfV_k$ respectively form $\bfR^{-1}$-orthogonal and $\bfQ$-orthogonal (c.f., \cref{eq:QV-orth}) bases for the Krylov subspaces as follows
\[\begin{aligned}
    \text{Span}\{\bfU_k\} &= \mathcal{K}_k(\bfA(\bfQ + \widehat{\bfD}^{-1})\bfA^{\top}\bfR^{-1}, \bfc) \\
    \text{Span}\{\bfV_k\} &= \mathcal{K}_k(\bfA^{\top}\bfR^{-1}\bfA(\bfQ + \widehat{\bfD}^{-1}),  \bfA^{\top}\bfR^{-1}\bfc)
\end{aligned}
\]
respectively. This means that when the diagonal weighting matrices are fixed, FGGK turns into a Krylov subspace method.

\subsection{Selecting regularization parameters}
\label{ss_regpar}
In this section, we describe various methods for selecting regularization parameters $\lambda$ and $\alpha$ for the projected problem \cref{eq:projected}. Notice that the solution to this problem is given by
$$ \bff_k(\lambda,\alpha) = (\bfM_k\t \bfM_k + \lambda^2 \bfI + \alpha^2{\bfR}_{W,k}\t{\bfR}_{W,k})^{-1}\bfM_k\t\beta_1\bfe_1.$$
We write $\bff_k(\lambda,\alpha)$ to denote the explicit dependence of $\bff_k$ on the regularization parameters $\lambda$ and $\alpha$.

Let $\bfC_k(\lambda,\alpha) = (\bfM_k\t\bfM_k + \lambda^2 \bfI + \alpha^2{\bfR}_{W,k}\t{\bfR}_{W,k})^{-1}\bfM_k\t$, and denote the projected residual 
$$\bfr^{\text{proj}}_k(\lambda, \alpha) = \bfM_k\bff_k(\lambda, \alpha) - \beta_1\bfe_1.$$ 
We use hybrid regularization techniques to estimate the regularization parameters in the projected space similar to
\cite{cho2020hybrid}. However, the key difference is that~\eqref{eq:map} is no longer a Tikhonov type problem, so the techniques do not apply directly. Therefore, the parameter selection techniques used on the projected problem we present below are heuristics inspired by~\cite{cho2020hybrid}.  

In order to provide a benchmark for comparing the parameter selection methods, we define the optimal parameters as  \begin{equation}
    (\lambda^{\text{proj}}, \alpha^{\text{proj}}) = \argmin_{\lambda, \alpha}\norm[2]{\bfs_k(\lambda,\alpha) - \bfs_{\text{true}}}^2,
\end{equation}
where $\bfs_{\text{true}}$ is the true solution that is not available in practice. As was mentioned, we use it merely to test the performance of the parameter selection methods.  
We can select parameters $\lambda, \alpha$ using the unbiased \textit{predictive risk estimation (UPRE) method} for the projection problem, where
\begin{equation}
    (\lambda^{\text{proj}}, \alpha^{\text{proj}}) = \argmin_{\lambda, \alpha} \frac{1}{k}\norm[2]{\bfr^{\text{proj}}_k(\lambda,\alpha)}^2 + \frac{2}{k}\text{tr}(\bfM_k\bfC_k(\lambda, \alpha)) - 1.
    \label{eq:upre}
\end{equation}
Notice that the noise level of the problem should be included in the definition of $\bfR.$
Another common approach is to use the \textit{Discrepancy Principle (DP)}, where parameters $\lambda, \alpha$ are selected such that
\begin{equation}
    (\lambda^{\text{proj}}, \alpha^{\text{proj}}) = \argmin_{\lambda, \alpha} \Big| \norm[2]{\bfr^{\text{proj}}_k(\lambda,\alpha)}^2  - m\tau \Big| ,
    \label{eq:dp}
\end{equation}
where 
$\tau \geq 1$ is a safety factor. 
Without a priori knowledge of the noise level, another option is to use an extension of the \textit{weighted generalized cross validation (WGCV)} method. The basic idea is to select parameters,
\begin{equation}
    (\lambda^{\text{proj}}, \alpha^{\text{proj}}) = \argmin_{\lambda, \alpha}\frac{\norm[2]{\bfr^{\text{proj}}_k(\lambda,\alpha)}^2}{(\text{tr}(\bfI_k - \omega \bfM_k\bfC_k(\lambda,\alpha)))^2},
    \label{eq:wgcv}
\end{equation}
where $\omega = k/m$ \cite{renaut2017hybrid}.
Various parameter choice methods can be used here since the projected problem is small.  We point the interested reader to other works on regularization parameter selection in the context of hybrid methods \cite{chung2019flexible}.

Motivated by the approaches described in \cite{chung2008weighted, chung2015hybrid}, we introduce three stopping criteria for the FGGK process in the solution decomposition hybrid approach. The iterative process is terminated if either of these conditions is satisfied: 
 (i) a maximum number of iterations is attained, (ii) the GCV function defined in terms of the iteration, 
\begin{equation}
\label{eq:GCVstop}
    \widehat{G}(k) = \frac{k\norm[2]{\bfr^{\text{proj}}_k(\lambda,\alpha)}^2}{(\text{tr}(\bfI_k -  \bfM_k\bfC_k(\lambda,\alpha)))^2}
\end{equation}
reaches the minimum or flattens out. In addition to these criteria, one could also consider stopping if the gradient of the objective function \eqref{eq:objfun} is sufficiently small; however, we did not need this in our implementation.

\section{Numerical results}
\label{sec:numerics}
In this section, we investigate the performance of the proposed \texttt{sdHybr} method using various examples from imaging processing. 
In \cref{sec:AT}, we consider a hypothetical atmospheric transport problem where the goal is to recover emissions maps that cover North America for detecting anomalies. Then, in \cref{sec:SDMT} we consider dynamic spherical means tomography reconstruction, where the true image combines moving smooth components and sparsely positioned dots. Finally, in \cref{sec:OCO2} we consider a more challenging case study from NASA's Orbiting Carbon Observatory 2 (OCO-2) satellite that includes an atmospheric transport model and dynamic data.  For the last example, the true solution is not synthetically generated as a sum of random fields.  For each of the case studies, we show that \texttt{sdHybr} methods are able to capture both smooth and sparse components.

We compare the performance of \texttt{sdHybr} methods to that of generalized hybrid methods and flexible hybrid algorithms, denoted by  \texttt{genHyBR} and \texttt{fHybr} respectively.
For the hybrid methods, we consider regularization parameters selected using the UPRE method, the DP method and the WGCV method.  For \texttt{sdHybr}, this corresponds to solving nonlinear constrained optimization problems \eqref{eq:upre}, \eqref{eq:dp} and \eqref{eq:wgcv} respectively. For this task, we use a Quasi-Newton method as implemented in MATLAB's \texttt{fminunc} function with an initial guess of $\lambda = -0.5$ and $\gamma = -0.5$. For the stopping criteria for \texttt{sdHybr}, the iterative method is terminated if either of the following two criteria is satisfied: (i) a maximum number of iterations is reached; (ii) the GCV stopping function $\widehat{G}$ defined in \eqref{eq:GCVstop} reaches the minimum or flattens out. The following experiments ran on a laptop
computer with Intel i5 CPU 2GHz and 16G memory.

\subsection{Case study 1: A hypothetical atmospheric transport problem}
\label{sec:AT}
\begin{figure}[ht!]
    \centering
    \includegraphics[width =\textwidth]{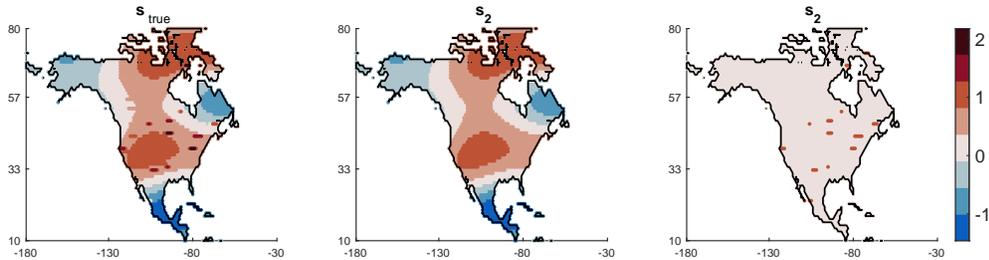}
    \caption{Atmospheric transport example, case study 1. The true emissions image $\bfs_{\rm true}$ provided on the left is the sum of smooth image $\bfs_1$ and sparse image $\bfs_2$, i.e., $\bfs_{\rm true} = \bfs_1 + \bfs_2$. Colormaps for all images are the same.
    }
    \label{fig:ex1_case1_true}
\end{figure}

We investigate a synthetic atmospheric transport problem, where observations $\bfd \in \bbR^{98880}$ are generated as in \eqref{eq:sdproblem} with $\bfA \in \bbR^{98880 \times 3222}$ representing a forward atmospheric model and $\bfs$ representing the true emissions map which is a summation of a randomly-generated smooth image and an image with sparse anomalies. Henceforth, we denote the true emissions vector as $\bfs_{\rm true}\in \bbR^{3222}$.  See \cref{fig:ex1_case1_true}. The smooth image is generated by a Mat$\Acute{\text{e}}$rn kernel ~\cite[Equation (4.14)]{williams2006gaussian} with parameters $\nu = 2.5$ and $\ell = 0.05$. For the image with sparse speckles,  a few of the speckles acquire the maximum value, whereas the remaining have varying values.  The goal is to reconstruct the unknown set of states or fluxes in space, where the spatial resolution is $1^{\text{o}} \times 1^{\text{o}}$.  This resolution is coarser than ideal for detecting super-emitters but provides a nice testbed example. In this case study, the observations $\bfd$ are sampled at the locations and times of OCO-2 observations during July through mid-August 2015, and the atmospheric model $\bfA$ is from NOAA's CarbonTracker-Lagrange project \cite{MillerSaibaba2020,liu2021}. Specifically, these atmospheric modeling simulations are from the Weather Research and Forecasting (WRF) Stochastic Time-Inverted Lagrangian Transport Model (STILT) modeling system \cite{Lin2003,Nehrkorn2010}. Note that we do not use realistic CO$_2$ emissions in this case study (c.f., case study 3).  Instead, we use randomly-generated emissions to create a relatively simple, initial test case for the algorithms proposed here. We add Gaussian white noise corresponding to a $4\%$ noise level to the observations, i.e., $\frac{\norm[2]{\bfdelta}}{\norm[2]{\bfA \bfs_{\text{true}}}} = 0.04$. 

\begin{figure}[ht!]
    \centering
    \includegraphics[width = \textwidth]{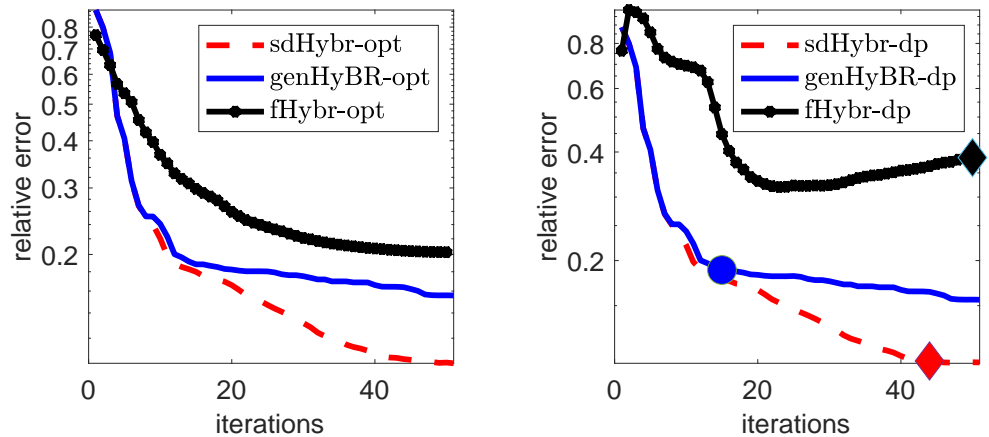}
    \caption{Atmospheric transport, case study 1: Relative reconstruction error norms per iteration of \texttt{sdHybr}, \texttt{genHyBR}, and \texttt{fHybr}.  Results in the left plot correspond to selecting the optimal regularization parameters at each iteration, and results in the right plot correspond to the DP-selected regularization parameters.}
    \label{fig:ex1_case1_err}
\end{figure}

We obtain reconstructions using the proposed \texttt{sdHybr} method and provide relative reconstruction error norms per iteration.  These are computed as \[\norm[2]{\bfs_k - \bfs_{\text{true}}}/\norm[2]{\bfs_{\text{true}}},\] where $\bfs_{k}$ is the reconstruction at the $k$th iteration. 
In the left plot of \cref{fig:ex1_case1_err}, we provide relative reconstruction error norms per iteration using the optimal regularization parameter at each iteration, which is not available in practice.  In the right plot, we see that similar results are obtained using the DP-selected regularization parameters. Results for \texttt{sdHybr} with WGCV and UPRE are very similar, so we do not provide them here. For comparison, we provide the relative reconstruction error norms per iteration of \texttt{genHyBR} and \texttt{fHybr} for both the optimal and DP-selected regularization parameters. For \texttt{genHyBR} and \texttt{sdHybr}, we let $\bfQ$ represent a Mat$\acute{\text{e}}$rn kernel with $\nu = 0.5$ and $\ell = 0.5$. All considered hybrid methods include a variety of stopping criteria. The tolerance for the GCV function was set to $\delta_{\text{GCV}}= 10^{-6}$ and the maximum number of iterations is $50$.  The diamonds denote the (automatically-selected) stopping iterations. \texttt{sdHybr} with different parameter selection resulted in similar stopping points.
 
\begin{figure}[ht!]
    \centering
    \includegraphics[width = \textwidth]{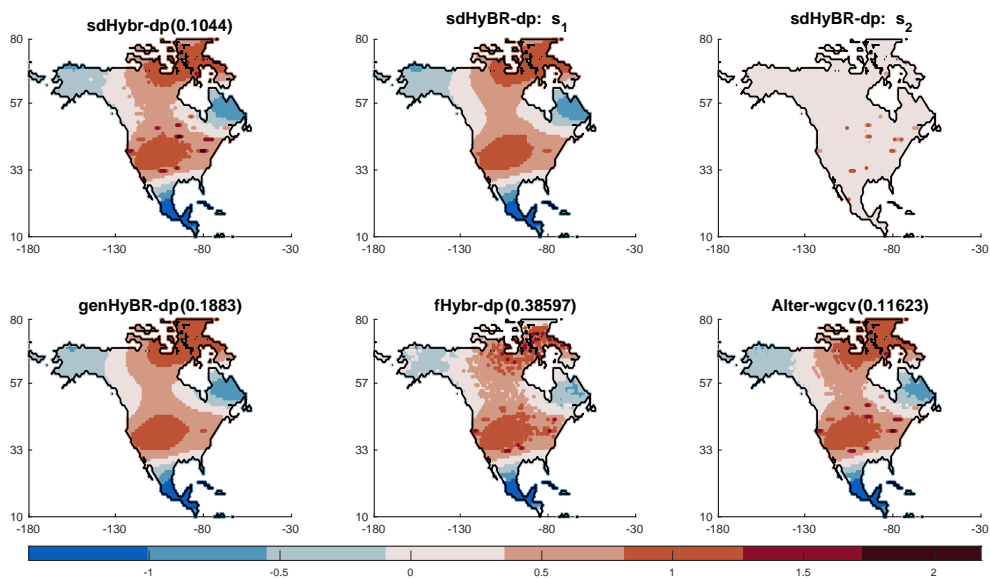}
    \caption{Atmospheric transport, case study 1: Reconstructions with relative reconstruction error norms provided in the titles. Decomposed solutions computed using \texttt{sdHybr-dp}.}
    \label{fig:ex1_case1_recon}
\end{figure}

We observe that \texttt{sdHybr} reconstructions produce smaller relative reconstruction errors than \texttt{genHyBR} and \texttt{fHybr}, demonstrating that our solution decomposition can be beneficial.  This result is also evident in the image reconstructions displayed in \cref{fig:ex1_case1_recon}.
In the top row of \cref{fig:ex1_case1_recon}, we provide the overall \texttt{sdHybr-dp} reconstruction, along with the computed estimates of $\bfs_1$ and $\bfs_2$.  In the bottom row, we provide the reconstructions obtained using \texttt{genHyBR-dp} and \texttt{fHybr-dp}.
We observe that the \texttt{genHyBR} reconstruction captures the smooth regions but fails to reconstruct the anomalies, while the \texttt{fHybr} reconstruction captures the anomalies but lacks smoothness in the background.

In the bottom row of \cref{fig:ex1_case1_recon}, we also provide a reconstruction using
\texttt{Alter-wgcv}, which is an alternating optimization approach \cref{alg:altopt}. Both \texttt{sdHybr}~and \texttt{Alter-wgcv} capture the sparse anomalies as well as the smooth background.
The maximum number of iterations for \texttt{Alter-wgcv} was $200$.
We found that the \texttt{Alter} approach requires high accuracy of the algorithms used to compute solutions in the alternating framework, hence the larger number of iterations for each solve within \texttt{Alter} and the overall longer CPU times: \texttt{sdHybr}~took $19.15$ seconds and \texttt{Alter} took $253.97$ seconds.  Furthermore, we observed that \texttt{Alter} is very sensitive to the accuracy of the initial guess.  Moreover, only WGCV was able to provide \texttt{Alter} reconstructions that could distinguish anomalies.

\subsection{Case study 2: Dynamic spherical means tomography}
\label{sec:SDMT}
In this experiment, we consider a dynamic tomography setup where the goal is to reconstruct a sequence of images from a sequence of projection datasets.  Such scenarios are common in dynamic photoacoustic or dynamic electrical impedance tomography, where the underlying parameters change during the data acquisition process \cite{wang2014fast,schmitt2002efficient2,hahn2014efficient}. Reconstruction is particularly challenging for nonlinear or nonparametric deformations and often requires including a spatiotemporal prior \cite{chung2018efficient,schmitt2002efficient1}.  
In spatiotemporal inversions, classic approaches (e.g., those based on parametric covariance families or separable covariance functions) may not be rich enough to capture the phenomena of interest alone, and multiple priors may be required to enforce different spatiotemporal properties. 

\begin{figure}[ht!]
    \centering
    \includegraphics[width = \textwidth]{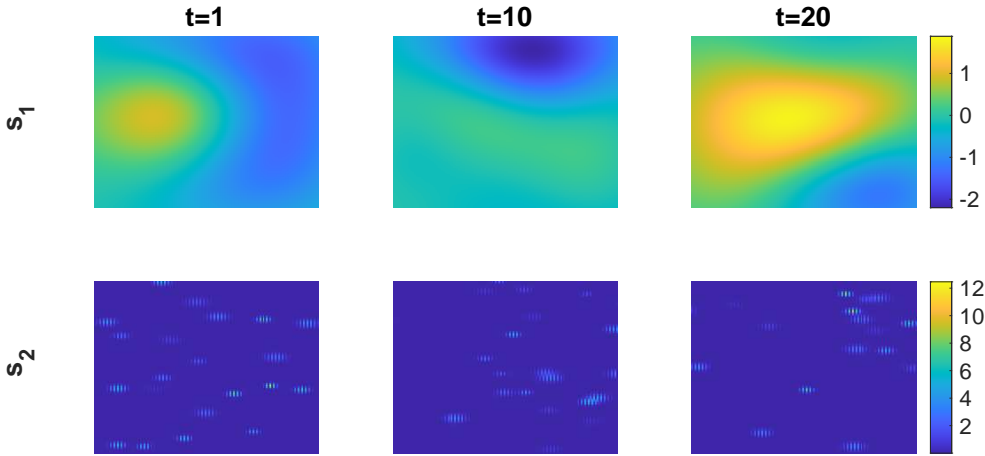}
    \caption{Dynamic spherical tomography example, case study 2. We provide three true images corresponding to time points $t=1, 10, 20.$ For each time point the true image is a sum of a true image plus a sparse image.}
    \label{fig:dynamic_true}
\end{figure}

In this example, we consider a sequence of $20$ true images (e.g., time points), where each image is $128\times 128$ and represents a sum of a smooth image and a sparse image.  That is, the true image at the $t$-th time point can be represented as the sum of two images: 
 the smooth image $\bfs_1$ was generated using a truncated Karhunen-Lo\'eve expansion using 30 basis vectors, with a Mat\'ern covariance kernel defined with two spatial and one temporal dimensions (see, for example,~\cite[Section]{chung2018efficient}). We will refer to this as a 3D Mat\'ern kernel. We also take $\nu = 0.2$ and $\ell = 0.2$. The sparse image $\bfs_2$ was generated using a star cluster example, and the two images were summed together.  
 In \cref{fig:dynamic_true}, we provide three of the true image decompositions (i.e., for time points $t=1, 10, 20$). Notice that although sparsely distributed, the spots in $\bfs_2$ take pixel values in a larger range compared to pixel values in $\bfs_1$. 

\begin{figure}[ht!]
    \centering
    \includegraphics[width = .75\textwidth]{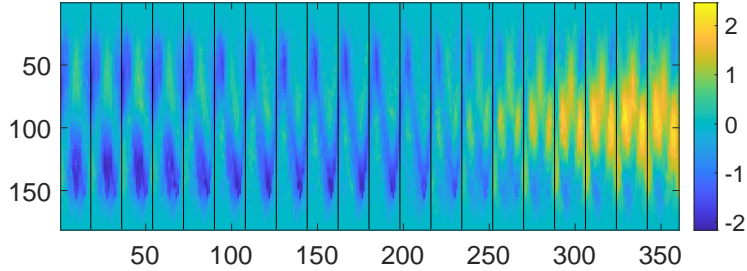}
    \caption{Dynamic spherical tomography example, case study 2. There are $20$ observed sinograms corresponding to $20$ time points.}
    \label{fig:dynamic_observe}
\end{figure}

We consider a linear problem of the form \cref{eq:sdproblem}, where
where $$\bfA =  \begin{bmatrix}
    \bfF_1 & & \\
 &\ddots &\\
 &&\bfF_{20}
\end{bmatrix}  \in \bbR^{{18*181} \times 128^2},$$ 
where $\bfF_t$ represents a spherical projection matrix corresponding to $18$ equally spaced angles between $t$ and $340+t$ for $t=1, \ldots, 20,$ and $\bfd \in \bbR^{20*18*181}$ contains projection data. For this, we use the \verb|PRspherical| test problem from the IRTools toolbox~\cite{gazzola2018ir,hansen2017air}, and to simulate measurement error we add $2\%$ Gaussian noise.  The collection of $20$ observed sinograms are concatenated and provided in \cref{fig:dynamic_observe}. For the reconstructions, we used prior covariance matrix $\bfQ$ representing a 3D Mat{\' e}rn kernel with $\nu = 0.5$ and $\ell = 0.4$.

\begin{figure}[ht!]
    \centering
    \includegraphics[width = \textwidth]{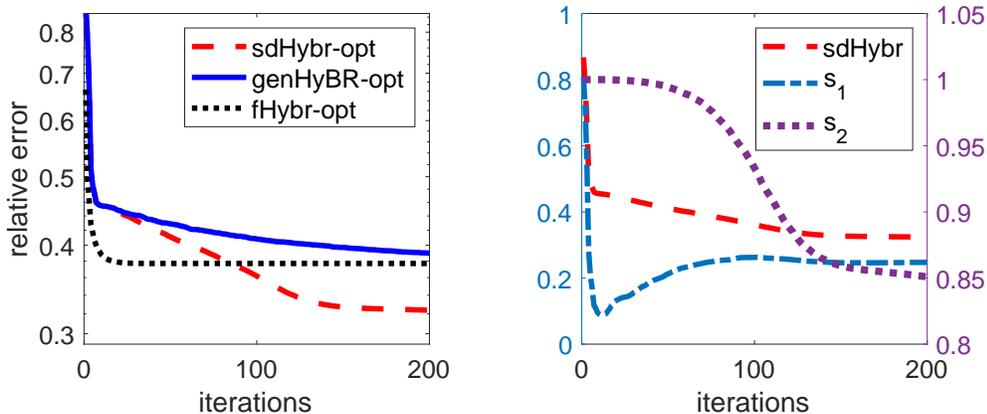}
    \caption{Dynamic spherical tomography example, case study 2. Relative reconstruction error norms per iteration are provided in the left plot for \texttt{sdHybr-opt}, \texttt{genHyBR-opt}, and \texttt{fHybr-opt}. In the right plot, we provide the relative reconstruction errors for $\bfs_1$ (left axis) and $\bfs_2$ (right axis) separately.  The relative reconstruction errors for \texttt{sdHybr-opt} from the left plot are provided again for reference.}
    \label{fig:dynamic_enrm}
\end{figure}

We focus on a comparison of methods using optimal regularization parameters, and we provide relative reconstruction error norms computed per iteration for \texttt{sdHybr-opt}, \texttt{genHyBR-opt}, and \texttt{fHybr-opt} in the left plot of \cref{fig:dynamic_enrm}. We observe that \texttt{sdHybr-opt} can achieve smaller overall reconstruction error norms compared to \texttt{genHyBR-opt} and \texttt{fHybr-opt}.  In the right plot of  \cref{fig:dynamic_enrm}, we provide relative reconstruction errors for $\bfs_1$ and $\bfs_2$ separately.  An interesting observation is that in early iterations, the \texttt{sdHybr-opt} method seems to reconstruct better approximations of $\bfs_1$, and in later iterations, reconstructions seem to capture features in $\bfs_2$.

\begin{figure}[ht!]
    \centering
    \includegraphics[width = \textwidth]{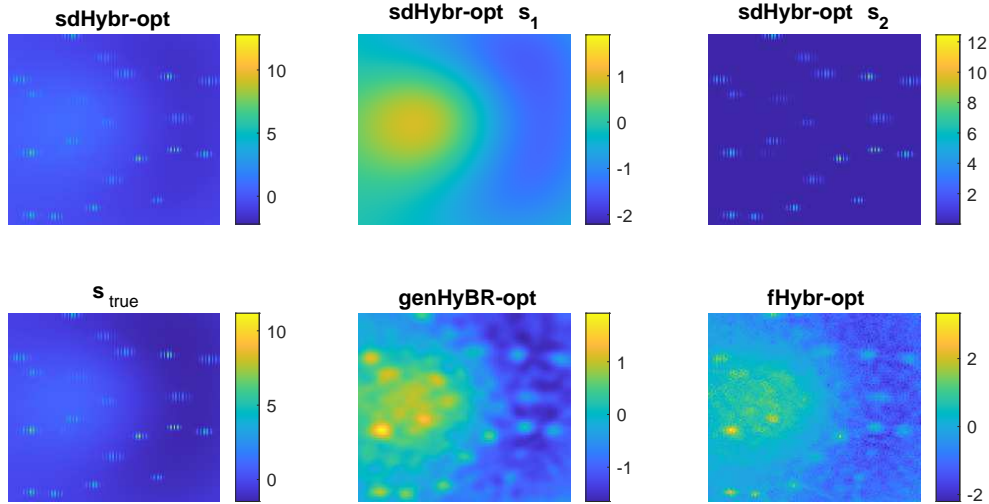}
    \caption{Dynamic spherical tomography example, case study 2. For time point $t=1$, the \texttt{sdHybr-opt} solution is provided along with reconstructed solutions $\bfs_1$ and $\bfs_2$ in the top row.  In the bottom row, the true image along with \texttt{genHyBR}, and \texttt{fHybr} reconstructions with the optimal regularization parameter are provided.}
    \label{fig:dynamic_recon}
\end{figure}

The main benefit of the solution decomposition approach can be  seen in the reconstructions. For time point $t=1$, we provide image reconstructions in \cref{fig:dynamic_recon}, along with the corresponding true image. 
We observe that \texttt{sdHybr-opt} can reconstruct better solutions than \texttt{genHyBR-opt} and \texttt{fHybr-opt}. Moreover, the solution decomposition approach can simultaneously solve for both components, which means that we have two separate image reconstructions $\bfs_1$ and $\bfs_2$. Notice that the \texttt{sdHybr-opt} solution can better capture the larger values in $\bfs_2$ while simultaneously capturing the smooth features in $\bfs_1$. 

We obtained similar results for the automatic parameter selection techniques, so we omit those here.  Also, we acknowledge that the choice of hyperparameters will be important for the overall reconstruction. 
Additional investigation is necessary to determine appropriate hyperparameters for the prior covariance matrix, and this is a topic of future work. 

\subsection{Case study 3: Atmospheric inverse modeling based on NASA's OCO-2 satellite}
\label{sec:OCO2}
For this case study, we consider a more realistic atmospheric inverse model, where the goal is to estimate CO$_2$ fluxes across North America using observations from NASA's OCO-2 satellite. The setup parallels that in \cite{MillerSaibaba2020,liu2021}, so we just provide an overview here.

We consider a linear model of the form \cref{eq:sdproblem}, where the aim is to estimate CO$_2$ fluxes at $3$-hourly temporal resolution over $41$ days (approximately $6$ weeks from late June through July 2015) and at $1^\circ\times1^\circ$ latitude-longitude spatial resolution.  This setup corresponds to $3,222$ unknowns per $3$-hour time interval; hence,
$\bfs \in \bbR^{328 \cdot 3222}$.
For $\bfs_{\rm true}$, we use CO$_2$ fluxes from NOAA's CarbonTracker product (version 2019b). 
Although a decomposition of $\bfs_{\rm true} = \bfs_1 + \bfs_2$ is not available, we observe that, similar to actual atmospheric models, the true fluxes contain a combination of large, sparsely distributed values which correspond to anomalies (e.g., fires, anthropogenic emissions, or anomalies in biospheric fluxes) and smooth, broad regions of surface fluxes with small-scale variability. Synthetic satellite observations given in $\bfd \in \bbR^{19,156}$ are generated as in \cref{eq:sdproblem}, where $\bfA$ represents the atmospheric transport simulation described in \cref{sec:AT} and $\bfdelta$ is added noise to represent measurement errors. The noise covariance matrix $\bfR$ is $\sigma^2 \bfI$, where $\sigma=0.5648$, which corresponds to a noise level of $50\%$. More specifically, for $\bfn\sim\calN(\bfzero,\bfI)$ the noise level of the observation corresponds to adding $\bfepsilon = \sigma\bfn$ where $\sigma = \texttt{nlevel}\cdot \frac{\|\bfA \bfs_{\rm true}\|_2}{\|\bfn\|_2}$.  Notice that this inverse problem is significantly under-determined, and thus, appropriate prior information plays a key role. 

Similar to previous approaches \cite{yadav2016statistical, MillerSaibaba2020}, we consider prior covariance matrix, $\bfQ = \lambda^{-2} \bfQ_t \kron \bfQ_s$ where $\bfQ_t$ represents the temporal covariance and $\bfQ_s$ represents the spatial covariance in the fluxes.  These covariance matrices are  defined by kernel functions \begin{align}
    k_t(d_t;\theta_t) & = \left\{ \begin{array}{ll}
         1-\frac{3}{2}\left(\frac{d_t}{\theta_t}\right)+\frac{1}{2}\left(\frac{d_t}{\theta_t}\right)^3  & \mbox{ if } d_t\leq \theta_t,\\
         0 & \mbox{ if } d_t > \theta_t,
    \end{array} \right. \\
    k_s(d_s;\theta_s) & = \left\{ \begin{array}{ll}
         1-\frac{3}{2}\left(\frac{d_s}{\theta_s}\right)+\frac{1}{2}\left(\frac{d_s}{\theta_s}\right)^3  & \mbox{ if } d_s\leq \theta_s,\\
         0 & \mbox{ if } d_s > \theta_s,
    \end{array} \right. 
\end{align} where $d_t$ is day difference between two unknowns, $d_s$ is spherical distance between two unknowns, and $\theta_t,\theta_g$ are kernel parameters. In this setting, we set $\theta_t=9.854$ and $\theta_s = 555.42$, as in \cite{MillerSaibaba2020}. 

\begin{figure}[bt]
    \centering
    \includegraphics[width = \textwidth]{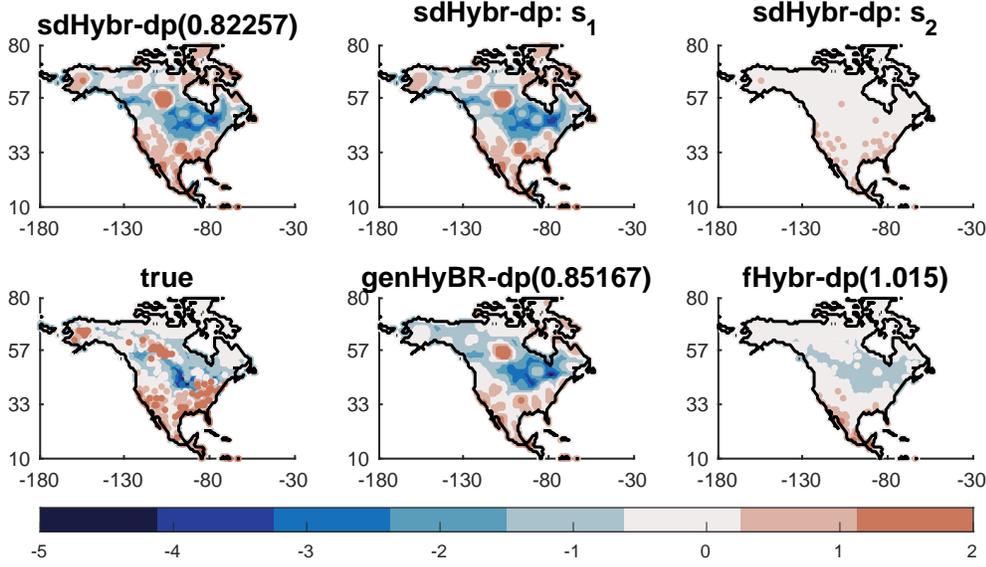}
    \caption{OCO-2 example, case study 3. In the top row, we provide the averaged computed reconstruction for \texttt{sdHybr}, along with the average reconstructions for $\bfs_1$ and $\bfs_2$.  In the bottom row, we provide the averaged true fluxes for reference, along with the reconstructions obtained using \texttt{genHyBR} and \texttt{fHybr}. Relative reconstruction error norms for the spatiotemporal fluxes are provided in the titles, and all results correspond to using the DP selected regularization parameters.}
    \label{fig:ex3_recon}
\end{figure}

We compute spatiotemporal reconstructions using \texttt{sdHybr}, \texttt{genHyBR}, and \texttt{fHybr}, all using the automatically selected regularization parameters using the discrepancy principle. We provide in \cref{fig:ex3_recon} the overall average image of flux reconstructions, along with the average true image. We observe that \texttt{genHyBR} does fairly well at estimating the broad regions of flux estimates, and \texttt{fHybr} is not able to capture a good reconstruction.  The average reconstruction of our proposed \texttt{sdHybr} method results best captures of both sources and sinks presented in the true average image.  Furthermore, a significant benefit of \texttt{sdHybr} is the ability to obtain the solution decomposition. The reconstructions of $\bfs_1$ and $\bfs_2$ from \texttt{sdHybr} are provided in the top row of \cref{fig:ex3_recon}.

The results of this case study demonstrate that \texttt{sdHybr} can yield accurate results for complex, spatiotemporal atmospheric inverse modeling.  The proposed method enables anomaly detection, due to the simultaneous reconstruction of two separate reconstructions, one containing smooth, broad regions and another capturing sparsely distributed anomalies. Moreover, the use of efficient hybrid projection methods means that these methods can be paired with the adjoint of an atmospheric model (where explicit construction of $\bfA$ is replaced with efficient matrix-vector products with $\bfA$ and $\bfA\t$) and with automatic selection of regularization parameters.

\section{Conclusions}
\label{sec:conclusions}
We have described hybrid projection methods for efficiently computing solutions to large-scale inverse problems, where the solution consists of two characteristically different solutions.  We focus on the scenario where the desired solution is a sum of a sparse solution and a smooth solution (e.g., a framework that can be used for anomaly detection) and describe a flexible, generalized Golub-Kahan hybrid iterative approach that confers several advantages. The approaches are \textit{efficient}, in part because they converge quickly, they exploit efficient matrix-vector multiplications, and they avoid expensive inner-outer optimization schemes by leveraging recent work on flexible preconditioning techniques. These methods are also \textit{automatic} since hyperparameters and stopping criteria can be determined as part of the iterative algorithm.  We describe various problems and alternative formulations that also fit within our framework, so these methods can be utilized for a wide range of problems. Numerical results from various applications, including dynamic inverse problems and atmospheric inverse modeling, demonstrate the benefits and potential for our approach.

Future work includes extending and applying the proposed method for real satellite data, that when coupled with the GEOS-Chem chemical transport model, requires a reformulation to handle the adjoint.  In addition, techniques for uncertainty quantification can be extended to this framework by utilizing a linearization approach to approximate the posterior at the MAP estimate with a Gaussian distribution  \cite{bardsley2018computational,bui2013computational}  or by exploiting previously computed bases for the solution subspaces as described in \cite{saibaba2020efficient}.

\appendix

\section{Alternative to FGGK}\label{subsec:app1}
   Here we briefly derive an alternative to FGGK. Recall that we have to solve the sequence of systems~\eqref{eq:MM_min} where now $\bfD_k = \bfD(\bfxi^{(k)})$. Let $\bfy = \bfD_k\bfxi$, so that the optimization problem can be recast as 
\begin{equation}
	\min_{\bfx,\ \bfy}  \norm[2]{\widehat \bfA \widehat \bfQ \begin{bmatrix} \bfI &  \\ & \bfD_k^{-1} \end{bmatrix}
	\begin{bmatrix} \bfx \\  \bfy \end{bmatrix} - \bfc}^2 +  \norm[\widehat \bfQ]{\begin{bmatrix} \lambda \bfI & \\ & \alpha \bfI \end{bmatrix} \begin{bmatrix} \bfx \\ \bfy \end{bmatrix}}^2
\end{equation}
where $\widehat \bfA$ and $\widehat{\bfQ}$ were defined in~\eqref{eqn:qhat}. We can derive a similar process as before, where after $k$ iterations of this process, we obtain 
\begin{equation}
	\widehat \bfA \widehat \bfQ \bfZ_k = \bfU_{k+1} \bfG_k \quad \mbox{and} \quad \widehat \bfA\t \bfU_{k+1} = \bfV_{k+1} \bfH_{k+1},
\end{equation}
where
\begin{itemize}
\item $\bfZ_k = \begin{bmatrix} \bfz_1 & \dots & \bfz_k \end{bmatrix} = \begin{bmatrix} \bfL_1  \bfv_1 & \dots & \bfL_k \bfv_k \end{bmatrix}\in \bbR^{2n \times k}$ with $\bfL_i = \begin{bmatrix}\bfI & \\ & \bfD^{-1}_i\end{bmatrix} \in \mathbb{R}^{2n \times 2n}$ contains solution basis vectors,
\item $\bfG_k =[g_{i,j}]_{{i=1,\dots,k+1; j=1,\dots,k}}\in \bbR^{(k+1)\times k}$ is upper Hessenberg,
\item $\bfH_{k+1}=[h_{i,j}]_{i,j=1,\dots,k+1}\in \bbR^{(k+1)\times (k+1)}$ is upper triangular,
\item $\bfU_{k+1} \!\!=\!\!
\begin{bmatrix}
 \bfu_1 & \dots & \bfu_{k+1}
\end{bmatrix}\!\in\!\bbR^{m \times (k+1)}$
with $\bfu_1=\bfc/\norm[2]{\bfc}$, and 
\item $\bfV_{k+1} =\begin{bmatrix}
 \bfv_1 & \dots & \bfv_{k+1}
\end{bmatrix}  \in \bbR^{2n\times (k+1)}$
such that $$\bfU_{k+1}\t \bfR^{-1}\bfU_{k+1} = \bfI_{k+1}\qquad \bfV_{k+1}\t \widehat\bfQ \bfV_{k+1} =\bfI_{k+1}.$$
\end{itemize}
The specific algorithm is given in \cref{alg:flexGK2}.

\begin{algorithm}[!ht]
\begin{algorithmic}[1]
\STATE Initialize $\bfu_1 = \bfc/\norm[\bfR^{-1}]{\bfc}$
\FOR {$i=1,\ldots,k$}
\STATE $\bfw = \widehat\bfA\t \bfR^{-1} \bfu_i$, $h_{j,i} = \bfw\t \widehat \bfQ \bfv_j$ for $j=1, \ldots, i-1$
\STATE $\bfw = \bfw - \sum_{j=1}^{i-1} h_{j,i}  \bfv_j$, $h_{i,i} = \norm[\widehat \bfQ]{\bfw}$, $\bfv_i = \bfw/h_{i,i}$
\STATE $\bfz_i = \begin{bmatrix}\bfI & \\ & \bfD^{-1}_i\end{bmatrix} \bfv_i$
\STATE $\bfw = \widehat \bfA \widehat \bfQ \bfz_i$, $g_{j,i} = \bfw\t \bfR^{-1}\bfu_j$ for $j=1, \ldots, i$
\STATE $\bfw = \bfw - \sum_{j=1}^i g_{j,i} \bfy_j$, $g_{i+1,i} = \norm[\bfR^{-1}]{\bfw}$, $\bfu_{i+1} = \bfw/g_{i+1,i}$
\ENDFOR
\end{algorithmic}
\caption{Alternative flexible, generalized Golub--Kahan process}
\label{alg:flexGK2}
\end{algorithm}

Next, assume that we seek solution $\begin{bmatrix} \bfx \\ \bfy \end{bmatrix} \in \calR(\bfZ_k)$, i.e., $\begin{bmatrix} \bfx\\ \bfy \end{bmatrix}  = \bfZ_k \bff$ for some $\bff \in \bbR^k$.  Then the projected problem,
\begin{equation}
	\min_{\begin{bmatrix} \bfx \\ \bfy \end{bmatrix} \in \calR(\bfZ_k)} \norm[2]{\widehat \bfA \widehat \bfQ
	\begin{bmatrix} \bfx \\ \bfy \end{bmatrix} - \bfc}^2 +  \norm[\widehat \bfQ]{\begin{bmatrix} \lambda \bfI & \\ & \alpha \bfI \end{bmatrix} \begin{bmatrix} \bfx \\ \bfy \end{bmatrix}}^2
\end{equation}
becomes
\[\min_{\bff \in \bbR^k} \norm[2]{\bfG_k \bff - \norm[\bfR^{-1}]{\bfc}\bfe_1}^2 +  \norm[\widehat \bfQ]{\begin{bmatrix} \lambda \bfI & \\ & \alpha \bfI \end{bmatrix} \bfZ_k \bff }^2.\]
The important point here is that, while the projected problem is in a smaller dimensional space, it is not obvious how to simultaneously estimate the projected parameters $\lambda$ and $\alpha$. One choice would be to assume $\alpha = \lambda \lambda_{\alpha}$, where $\lambda_\alpha$ is fixed and estimate $\lambda$ using techniques similar to \cref{ss_regpar}. However, the choice we made in \cref{ss_hybridmethod} allows for the estimation of both parameters within the projected space, while maintaining the same computational cost. This highlights the novelty of the proposed approach.

\bibliographystyle{abbrv}
\bibliography{references}

\begin{thebibliography}{10}

\bibitem{ambikasaran2012large}
S.~Ambikasaran, J.~Li, P.~Kitanidis, and E.~Darve.
\newblock Large-scale stochastic linear inversion using {H}ierarchical
  matrices.
\newblock {\em Computational Geosciences}, 2012.

\bibitem{bardsley2018computational}
J.~M. Bardsley.
\newblock {\em Computational Uncertainty Quantification for Inverse Problems},
  volume~19.
\newblock SIAM, 2018.

\bibitem{beck2009fast}
A.~Beck and M.~Teboulle.
\newblock A fast iterative shrinkage-thresholding algorithm for linear inverse
  problems.
\newblock {\em SIAM J. Imaging Sci.}, 2(1):183--202, 2009.

\bibitem{bui2013computational}
T.~Bui-Thanh, O.~Ghattas, J.~Martin, and G.~Stadler.
\newblock A computational framework for infinite-dimensional {B}ayesian inverse
  problems {P}art {I}: The linearized case, with application to global seismic
  inversion.
\newblock {\em SIAM Journal on Scientific Computing}, 35(6):A2494--A2523, 2013.

\bibitem{chen2016elastic}
D.-H. Chen, B.~Hofmann, and J.~Zou.
\newblock Elastic-net regularization versus $\ell$-1-regularization for linear
  inverse problems with quasi-sparse solutions.
\newblock {\em Inverse Problems}, 33(1):015004, 2016.

\bibitem{cho2020hybrid}
T.~Cho, J.~Chung, and J.~Jiang.
\newblock Hybrid projection methods for large-scale inverse problems with mixed
  {G}aussian priors.
\newblock {\em Inverse Problems}, 2020.

\bibitem{chung2019flexible}
J.~Chung and S.~Gazzola.
\newblock Flexible {K}rylov methods for $\ell_p$ regularization.
\newblock {\em SIAM Journal on Scientific Computing}, 41(5):S149--S171, 2019.

\bibitem{chung2008weighted}
J.~Chung, J.~G. Nagy, and D.~P. O’Leary.
\newblock A weighted {GCV} method for {L}anczos hybrid regularization.
\newblock {\em Electronic Transactions on Numerical Analysis}, 28:149--167,
  2008.

\bibitem{chung2015hybrid}
J.~Chung and K.~Palmer.
\newblock A hybrid {LSMR} algorithm for large-scale {T}ikhonov regularization.
\newblock {\em SIAM Journal on Scientific Computing}, 37(5):S562--S580, 2015.

\bibitem{chung2017generalized}
J.~Chung and A.~K. Saibaba.
\newblock Generalized hybrid iterative methods for large-scale {B}ayesian
  inverse problems.
\newblock {\em SIAM Journal on Scientific Computing}, 39(5):S24--S46, 2017.

\bibitem{chung2018efficient}
J.~Chung, A.~K. Saibaba, M.~Brown, and E.~Westman.
\newblock Efficient generalized {G}olub--{K}ahan based methods for dynamic
  inverse problems.
\newblock {\em Inverse Problems}, 34(2):024005, 2018.

\bibitem{fornasier2016cg}
M.~Fornasier, S.~Peter, H.~Rauhut, and S.~Worm.
\newblock Conjugate gradient acceleration of iteratively re-weighted least
  squares methods.
\newblock {\em Computational optimization and applications}, 65(1):205--259,
  2016.

\bibitem{gazzola2018ir}
S.~Gazzola, P.~C. Hansen, and J.~G. Nagy.
\newblock {IR} {T}ools: a {MATLAB} package of iterative regularization methods
  and large-scale test problems.
\newblock {\em Numerical Algorithms}, pages 1--39, 2018.

\bibitem{gazzola2014generalized}
S.~Gazzola and J.~G. Nagy.
\newblock Generalized {A}rnoldi--{T}ikhonov method for sparse reconstruction.
\newblock {\em SIAM Journal on Scientific Computing}, 36(2):B225--B247, 2014.

\bibitem{IRNekki}
I.~F. Gorodnitsky and B.~D. Rao.
\newblock A new iterative weighted norm minimization algorithm and its
  applications.
\newblock In {\em IEEE Sixth SP Workshop on Statistical Signal and Array
  Processing}, pages 412--415, 1992.

\bibitem{hahn2014efficient}
B.~N. Hahn.
\newblock Efficient algorithms for linear dynamic inverse problems with known
  motion.
\newblock {\em Inverse Problems}, 30(3):035008, 2014.

\bibitem{hansen2017air}
P.~C. Hansen and J.~S. J{\o}rgensen.
\newblock {AIR} tools {II}: algebraic iterative reconstruction methods,
  improved implementation.
\newblock {\em Numerical Algorithms}, pages 1--31, 2017.

\bibitem{huang2017some}
Y.~Huang and Z.~Jia.
\newblock Some results on the regularization of {LSQR} for large-scale discrete
  ill-posed problems.
\newblock {\em Science China Mathematics}, 60(4):701--718, Apr 2017.

\bibitem{jin2009elastic}
B.~Jin, D.~A. Lorenz, and S.~Schiffler.
\newblock Elastic-net regularization: error estimates and active set methods.
\newblock {\em Inverse Problems}, 25(11):115022, 2009.

\bibitem{kaipio2006statistical}
J.~Kaipio and E.~Somersalo.
\newblock {\em Statistical and Computational Inverse Problems}, volume 160.
\newblock Springer Science \& Business Media, 2006.

\bibitem{kitanidis1986parameter}
P.~K. Kitanidis.
\newblock Parameter uncertainty in estimation of spatial functions: Bayesian
  analysis.
\newblock {\em Water Resources Research}, 22(4):499--507, 1986.

\bibitem{kitanidis1995quasi}
P.~K. Kitanidis.
\newblock Quasi-linear geostatistical theory for inversing.
\newblock {\em Water Resources Research}, 31(10):2411--2419, 1995.

\bibitem{kitanidis1983geostatistical}
P.~K. Kitanidis and E.~G. VoMvoris.
\newblock A geostatistical approach to the inverse problem in groundwater
  modeling (steady state) and one-dimensional simulations.
\newblock {\em Water resources research}, 19(3):677--690, 1983.

\bibitem{lange2016mm}
K.~Lange.
\newblock {\em MM optimization algorithms}.
\newblock SIAM, 2016.

\bibitem{li2016novel}
Z.~Li, A.~Desolneux, S.~Muller, and A.-K. Carton.
\newblock A novel {3D} stochastic solid breast texture model for x-ray breast
  imaging.
\newblock In {\em International Workshop on Breast Imaging}, pages 660--667.
  Springer, 2016.

\bibitem{Lin2003}
J.~C. Lin, C.~Gerbig, S.~C. Wofsy, A.~E. Andrews, B.~C. Daube, K.~J. Davis, and
  C.~A. Grainger.
\newblock A near-field tool for simulating the upstream influence of
  atmospheric observations: The stochastic time-inverted {L}agrangian transport
  ({STILT}) model.
\newblock {\em Journal of Geophysical Research: Atmospheres}, 108(D16), 2003.

\bibitem{liu2021}
X.~Liu, A.~L. Weinbren, H.~Chang, J.~Tadi\'c, M.~E. Mountain, M.~E. Trudeau,
  A.~E. Andrews, Z.~Chen, and S.~M. Miller.
\newblock Data reduction for inverse modeling: an adaptive approach v1.0.
\newblock {\em Geoscientific Model Development Discussions}, 2020:1--21, 2020.

\bibitem{matheron_1973}
G.~Matheron.
\newblock The intrinsic random functions and their applications.
\newblock {\em Advances in Applied Probability}, 5(3):439–468, 1973.

\bibitem{matsui2011variable}
H.~Matsui and S.~Konishi.
\newblock Variable selection for functional regression models via the l1
  regularization.
\newblock {\em Computational Statistics \& Data Analysis}, 55(12):3304--3310,
  2011.

\bibitem{MillerSaibaba2020}
S.~M. Miller, A.~K. Saibaba, M.~E. Trudeau, M.~E. Mountain, and A.~E. Andrews.
\newblock Geostatistical inverse modeling with very large datasets: an example
  from the orbiting carbon observatory 2 ({OCO}-2) satellite.
\newblock {\em Geoscientific Model Development}, 13(3):1771--1785, 2020.

\bibitem{Nehrkorn2010}
T.~Nehrkorn, J.~Eluszkiewicz, S.~Wofsy, J.~Lin, C.~Gerbig, M.~Longo, and
  S.~Freitas.
\newblock Coupled weather research and forecasting–stochastic time-inverted
  {L}agrangian transport ({WRF–STILT}) model.
\newblock {\em Meteorol Atmos Phys}, 107:51--64, 2010.

\bibitem{renaut2017hybrid}
R.~A. Renaut, S.~Vatankhah, and V.~E. Ardestani.
\newblock Hybrid and iteratively reweighted regularization by unbiased
  predictive risk and weighted {GCV} for projected systems.
\newblock {\em SIAM J. Sci. Comput.}, 39(2):B221--B243, 2017.

\bibitem{wohlberg2008lp}
P.~Rodriguez and B.~Wohlberg.
\newblock An efficient algorithm for sparse representations with $\ell^p$ data
  fidelity term.
\newblock In {\em Proceedings of 4th IEEE Andean Technical Conference
  (ANDESCON)}, 2008.

\bibitem{saibaba2020efficient}
A.~K. Saibaba, J.~Chung, and K.~Petroske.
\newblock Efficient {K}rylov subspace methods for uncertainty quantification in
  large {B}ayesian linear inverse problems.
\newblock {\em Numerical Linear Algebra with Applications}, 27(5):e2325, 2020.

\bibitem{schmitt2002efficient1}
U.~Schmitt and A.~K. Louis.
\newblock Efficient algorithms for the regularization of dynamic inverse
  problems: {I. Theory}.
\newblock {\em Inverse Problems}, 18(3):645, 2002.

\bibitem{schmitt2002efficient2}
U.~Schmitt, A.~K. Louis, C.~H. Wolters, and M.~Vauhkonen.
\newblock Efficient algorithms for the regularization of dynamic inverse
  problems: {II. Applications}.
\newblock {\em Inverse Problems}, 18(3):659, 2002.

\bibitem{tremoulheac2014dynamic}
B.~Tr{\'e}moulh{\'e}ac, N.~Dikaios, D.~Atkinson, and S.~R. Arridge.
\newblock Dynamic {MR} image reconstruction--separation from undersampled
  ($k,t$)-space via low-rank plus sparse prior.
\newblock {\em IEEE Transactions on Medical Imaging}, 33(8):1689--1701, 2014.

\bibitem{wang2014fast}
K.~Wang, J.~Xia, C.~Li, L.~V. Wang, and M.~A. Anastasio.
\newblock Fast spatiotemporal image reconstruction based on low-rank matrix
  estimation for dynamic photoacoustic computed tomography.
\newblock {\em Journal of Biomedical Optics}, 19(5):056007--056007, 2014.

\bibitem{williams2006gaussian}
C.~K. Williams and C.~E. Rasmussen.
\newblock {\em Gaussian processes for machine learning}, volume~2.
\newblock MIT press Cambridge, MA, 2006.

\bibitem{xiong2019convex}
K.~Xiong, G.~Zhao, G.~Shi, and Y.~Wang.
\newblock A convex optimization algorithm for compressed sensing in a complex
  domain: The complex-valued split {B}regman method.
\newblock {\em Sensors}, 19(20):4540, 2019.

\bibitem{yadav2016statistical}
V.~Yadav, A.~M. Michalak, J.~Ray, and Y.~P. Shiga.
\newblock A statistical approach for isolating fossil fuel emissions in
  atmospheric inverse problems.
\newblock {\em Journal of Geophysical Research: Atmospheres}, 121(20):12--490,
  2016.

\bibitem{yin2015minimization}
P.~Yin, Y.~Lou, Q.~He, and J.~Xin.
\newblock Minimization of $\ell_{1-2}$ for compressed sensing.
\newblock {\em SIAM Journal on Scientific Computing}, 37(1):A536--A563, 2015.

\bibitem{zou2005regularization}
H.~Zou and T.~Hastie.
\newblock Regularization and variable selection via the elastic net.
\newblock {\em Journal of the royal statistical society: series B (statistical
  methodology)}, 67(2):301--320, 2005.

\end{thebibliography}

\end{document}